\documentclass[fleqn]{article}
\usepackage{amsmath}
\usepackage{bm}
\usepackage{amssymb}
\usepackage{indentfirst}
\usepackage{mathrsfs}
\usepackage{setspace}
\usepackage{geometry}
\usepackage{tikz}
\RequirePackage{amsmath,amsthm,amssymb}
\theoremstyle{plain}
\providecommand{\theoremname}{Theorem}%

\providecommand{\axiomname}{Axiom}%

\providecommand{\lemmaname}{Lemma}%

\providecommand{\corollaryname}{Corollary}%

\providecommand{\assertionname}{Assertion}%

\providecommand{\propositionname}{Proposition}%

\providecommand{\conjecturename}{Conjecture}%

\theoremstyle{definition}
\providecommand{\definitionname}{Definition}%

\providecommand{\examplename}{Example}%

\theoremstyle{remark}
\providecommand{\remarkname}{Remark}%
\newtheorem*{remark}{\remarkname}
\newcommand{\enabstractname}{Abstract}
\newenvironment{enabstract}{%
    \par\small
    \noindent\mbox{}\hfill{\bfseries \enabstractname}\hfill\mbox{}\par
    \vskip 2.5ex}{\par\vskip 2.5ex}
     \author{Yannan Qiu\thanks{Partially supported by NSFC grant No.11621061.}}
     \title{Some translations in the second lowest two-sided cell \\ of an affine Weyl group}
\date{}
\begin{document}
\begin{spacing}{1.5}
\maketitle
{\footnotesize 
\centerline{School of Mathematical Sciences,  Zhejiang University}
 \centerline{Zhejiang 310058,  China}
 \centerline{Academy of Mathematics and Systems Science, Chinese Academy of Sciences}
 \centerline{Beijing 100190,  China}
 \centerline{qiuyannan@zju.edu.cn}}

\begin{enabstract}
We are interested in cell properties of translations in affine Weyl groups. We find out some translations in the second lowest two-sided cell of an affine Weyl group and use the translations to formulate a refinement of Jianyi Shi's conjecture on the number of left cells in the second lowest two-sided cell. We verify the refinement for an affine Weyl group of type $\tilde A_{n-1}$ and type $\tilde G_2$. 

  \end{enabstract}

\section{Introduction}
Kazhdan-Lusztig cells play an important role in Kazhdan-Lusztig theory. In the theory an  interesting topic is the study of cells in affine Weyl groups. Cell decomposition in affine Weyl groups has drawn much attention. But the picture is far more from clear. By now we know the cell decomposition for the following types : rank two of any type (Lusztig \cite{lusztig1985cells}), type $\tilde A_n$ (Shi \cite{shi1986kazhdan}), type $\tilde C_3$ (Bedard \cite{bedard1986cells}),  type $\tilde B_3$ (Du \cite{jie1988decomposition}), type $\tilde C_4, \tilde D_4, \tilde F_4$ (Shi \cite{shi1994left}, \cite{shi1998left}, \cite{jian1998left}), type $\tilde B_4$ (Zhang \cite{zhang1994cell}), type $\tilde D_5$ (Chen \cite{chengdong2001left}). For affine Weyl groups of other types, not much  is known about cell decomposition, but the second highest two-sided cell and the lowest two-sided cell of an affine Weyl group are well understood (see \cite{lusztig1983some}, \cite{shi1987two}, \cite{jian1988two}). Also in \cite{shi2011second} Shi studied the second lowest two-sided cell of an affine Weyl group and gave a conjecture on the number of left cells in the second lowest two-sided cell.

All translations in an affine Weyl group form a normal subgroup of the affine Weyl group and the normal subgroup has finite index in the affine Weyl group. This suggests that cell properties of translations in affine Weyl groups should be interesting. In \cite{nanhua1989approach}, Nanhua Xi obtained some cell properties of the translations and used the properties to prove that each left cell of an affine Weyl group has only finitely many connected components.

We are interested in cell properties of translations in affine Weyl groups. A specific question is to find out translations in a given two-sided cell and in a given left cell of an affine Weyl groups.  In this paper we find out some translations in the second lowest two-sided cell of an affine Weyl group (Theorems 4.5, 4.8, 4.12, 4.14, 4.17)  and use the translations to formulate a refinement of Jianyi Shi's conjecture on the number of left cells in the second lowest two-sided cell (Conjecture 5.2). We verify the refinement for an affine Weyl group of type $\tilde A_{n-1}$ or type $\tilde G_2$ (Theorems 5.7, 5.6). For type $\tilde A_{n-1}$, we can say more. In addition, we determine all translations in the lowest two-sided cell of an affine Weyl group (Theorem 3.8). 

\section{Affine Weyl group and extended affine Weyl group}

In this section we fix some notations and recall some known facts.

{\bf 1.1.} Let $E$ be an euclidean space of dimension $n$ and $\Phi\subset E$ be an irreducible root system which spans $E$. Fix a simple root system $\Delta=\{\alpha_1,\ \alpha_2,\ \cdots, \alpha_n\}$ of $\Phi$. Let $\Phi^+$ and $\Phi^-$ be the set of  positive roots and the set of negative roots  in $\Phi$ respectively. 

For any $\alpha\in \Phi$,  the corresponding reflection $s_\alpha$ on $E$ is defined by $s_\alpha(x)=x -\langle x, \alpha^\vee\rangle\alpha$, where $\alpha^\vee=2\alpha/\langle\alpha,\ \alpha\rangle$, and $\langle\cdot,\cdot\rangle$ is the inner product on $E$.  The reflections $s_\alpha\ (\alpha\in\Phi)$ generate the Weyl group $W_0$   of $\Phi$.

Let $Q$ denote the root lattice $\mathbb{Z}\Phi$. Let $P=\{\lambda\in E\,|\, \langle \lambda, \alpha^\vee\rangle\in \mathbb Z\text{\ for all}\ \alpha\in\Phi\}$ be the weight lattice. Then $Q$ is a subgroup of $P$ with finite index. The weight lattice  $P$ is a free abelian group. Define $\lambda_1,\ \lambda_2,\ ...,\ \lambda_n$ by the conditions $\langle \lambda_i, {\alpha_j}^\vee\rangle=\delta_{ij}$ for $i, j=1,2,...,n$. The weights $\lambda_1,\ \lambda_2,\ ...,\ \lambda_n$ are called the fundamental dominant weights (with respect to $\Delta$). The fundamental dominant weights form a basis of $P$.

{\bf 1.2.} For each element $\lambda$ in $P$, define the translation $t_\lambda: E\to E,\ x\mapsto x+\lambda$. Then the translations $t_\lambda,\ \lambda\in P$ form a transformation group $X$ of $E$. Let $N$ be the subgroup of $X$ consisting of all $t_\lambda,\ \lambda\in Q$. Clearly, the map $\lambda\mapsto t_\lambda$ defines a group isomorphism from $P$ to $X$.

Let $W_a$ be the transformation group of $E$ generated by  $W_0$ and $N$,  and let $W$ be the transformation group of $E$ generated by $W_0$ and $X$. Then $W_a$ is an affine Weyl group and $W$ is an extended affine group. We have $W_a=W_0\ltimes N$ and $W=W_0\ltimes X$. In particular, $W_a$ has finite index in $W$.

Let $-\alpha_0$ be the highest short root of $\Phi$. Define $s_0=t_{-\alpha_0}s_{\alpha_0}$ and $s_i=s_{\alpha_i},\ 1\leq i\leq n$. Let $S=\{s_0,\ s_1,\ ...,\ s_n\}$. Then $(W_a,S)$ and $(W_0,S_0)$ are   Coxeter systems, where $S_0=\{s_1,\ s_2,\ \cdots, s_n\}$. Let $l:W_a\to\mathbb N$ be the length function of $W_a$. There is a finite subgroup $\Gamma$ of $W$ such that $W=\Gamma\ltimes W_a$. Following Lusztig we extend the length function to $W$ by setting $l(\tau w)=l(w)$ for any $\tau\in\Gamma$ and $w\in W_a$. Similarly, we extend the Bruhat order $\le$ on $W_a$ to $W$ by setting $\tau w\le \tau' u$ if $\tau=\tau'$ and $w\le u$.  We denote the extensions again by $l$ and $\leq$ respectively.

{\bf 1.3.} Iwahori and Matsumoto established a nice formula for the length function of $W$ (see Proposition 1.23 in \cite{iwahori1965some}):
\begin{equation}
l(t_\lambda w)=\sum_{\substack{\alpha\in{\Phi^+}\\ w^{-1}(\alpha)\in{\Phi^-}}}{|\langle \lambda, \alpha^\vee\rangle -1|}+
\sum_{\substack{\alpha\in{\Phi^+}\\ w^{-1}(\alpha)\in{\Phi^+}}}{|\langle \lambda, \alpha^\vee\rangle |},\quad \text{for all } t_\lambda\in X,\ \ w\in W_0.
\end{equation}
The above formula is equivalent to the following formula:
\begin{equation}
l(wt_\lambda )=\sum_{\substack{\alpha\in{\Phi^+}\\ w(\alpha)\in{\Phi^-}}}{|\langle \lambda, \alpha^\vee\rangle +1|}+
\sum_{\substack{\alpha\in{\Phi^+}\\ w(\alpha)\in{\Phi^+}}}{|\langle \lambda, \alpha^\vee\rangle |},\qquad \text{for all } t_\lambda\in X,\ \ w\in W_0.
\end{equation}

The subgroup $\Gamma$ actually consists of elements in $W$ with length 0.  Define $x_i=t_{\lambda_i}$. Then $x_1,\ ...,\ x_n$ form a basis of the abelian group $X$. By abuse of the term of fundamental dominant weights, {\it we also call $x_i$ a fundamental dominant weight in $X$}. Using Iwahori-Matsumoto's formula (1) and (2), we see
\begin{alignat}{2} &l(s_ix_j)=l(x_js_i)=l(x_j)+1,\quad if\ 1\le i\ne j\le n,\\
&l(s_ix_i)=l(x_i)+1,\quad l(x_is_i)=l(x_i)-1,\quad if\ 1\le i \le n.\end{alignat}

Let $P^+=\{\lambda\in P\,|\, \langle\lambda,\alpha_i^\vee\rangle\ge 0\text{\ \ for \ } i=1,2,...,n\}$ be the set of dominant weights. Correspondingly set $X^+=\{t_\lambda\,|\,\lambda\in P^+\}$. {\it The elements in $X^+$ will also be called dominant weights}. Using Iwahori-Matsumoto's formula (2) or (1) we see the following three identities.
\begin{alignat}{3}&X^+=\{ x\in X\,|\, l(wx)=l(w)+l(x),\ \forall w\in W_0\}=\{x_1^{a_1}x_2^{a_2}\cdots x_n^{a_n}\,|\, a_1,a_2,...,a_n\in\mathbb N\}.\\
& l(xy)=l(x)+l(y)\quad \text{for all }\ x,\ y\in X^+;\\
& l(wx^kw^{-1})=kl(x)\quad \text{for all } \ x\in X,\ w\in W_0,
 \end{alignat}

\section{Kazhdan-Lusztig Cells and Lusztig's $a$-function}
In \cite{kazhdan1979representations} , Kazhdan and Lusztig introduced the concepts of left, right, two-sided cell for Coxeter groups. In \cite{lusztig1985cells},  Lusztig defined the $a$-function $a:W\rightarrow \mathbb{N}$, which is a useful tool to study cells in Coxeter groups. The concepts of left, right, two-sided cell for Coxeter groups and the $a$-function work well for extended affine Weyl groups (see \cite{Lusztig1987CellsIA} ).

From now on, we keep the notations in section 1. In particular, $W_0$, $W_a$, $W$ are the Weyl group, affine Weyl group, the extended affine Weyl group associated to a root system $\Phi$, respecively. Of course,  some concepts and facts in this section are valid for more general (extended) Coxeter groups, but those are not used in this paper.

{\bf 2.1.}  The preorders $\underset{L}\leq,\ \underset{R}\leq,\ \underset{LR}\leq$ and the associated equivalence relations $\underset{L}\sim,\ \underset{R}\sim,\ \underset{LR}\sim$ on $W_a$ are defined  in \cite{kazhdan1979representations}. The preorders  are extended to $W$ by Lusztig in \cite{Lusztig1987CellsIA} as follows, for $\tau,\tau'\in \Gamma,\ w,u\in W_a$, set
$$\tau w\underset{L}\leq \tau'u\Longleftrightarrow \tau=\tau',\  w\underset{L}\leq u ,$$
$$\tau w\underset{R}\leq \tau'u\Longleftrightarrow \tau=\tau',\  w\underset{R}\leq u ,$$
$$\tau w\underset{LR}\leq \tau'u\Longleftrightarrow \tau=\tau',\  w\underset{LR}\leq u .$$
Then the equivalence relations $\underset{L}\sim,\ \underset{R}\sim,\ \underset{LR}\sim$ on $W_a$ are extended to $W$.

Similarly, the $a$-function is extended to $W$ by setting $a(\tau w)=a(w)$ for any $\tau\in\Gamma,\ w\in W_a$. Actually Lusztig defined the preorders and $a$-function on $W$ using Hecke algebra of $W$.

The equivalence classes with respect to $\underset{L}\sim$ ( respectively $\underset{R}\sim, \ \underset{R}\sim$) are called left cells ( respectively right cells, two-sided cells). And the preorders on $W$ induce  partial orders on the set of left cells (respectively right cells, two-sided cells) of $W$.

{\bf 2.2.} For any $w\in W$, we associate two subsets of $S$ as follows:
$$\mathscr{L}(w)=\{s\in S\ |\ sw\le w\},\qquad\mathscr{R}(w)=\{s\in S\ |\ ws \le w\}.$$

 According to  \cite{kazhdan1979representations}, we have

(a) If  $w\underset{L}\leq u$ (respectively,  $w\underset{R}\leq u$), then $\mathscr{R}(w)\supset \mathscr{R}(u)$ (respectively,  $\mathscr{L}(w)\supset \mathscr{L}(u)$). In particular, if  $w\underset{L}\sim u$( respectively,  $w\underset{R}\sim u$), then $\mathscr{R}(w)= \mathscr{R}(u)$( respectively,  $\mathscr{L}(w)=\mathscr{L}(u)$).

\medskip

For convenience, we often write $w\cdot u$ for $wu$ if $l(wu)=l(w)+l(u)$. Recall that $x_i=t_{\lambda_i}$.

{\bf Lemma 2.3.} Let $J\subset \{1,\ 2,\ \cdots,\ n\}$ and $x =\prod_{i\in J}x_i^{a_i},\ a_i\ge 1$, then

(a) $\mathscr{R}(x)=\{s_i\,|\, i\in J\}$;

(b) $x=(xw_J)\cdot w_J$, i.e. $l(x)=l(xw_J)+l(w_J)$, where $w_J$ is the longest element of the parabolic subgroup  $W_J=\langle s_i\,|\, i\in J\rangle$ of $W$.

\begin{proof}
The first result follows from formula (4) and (6), while the second result follows from the length formula (1) of Iwahori-Matsumoto,
$$l(xw_J)=\sum_{\substack{\alpha\in {\Phi}^+\\ {w_J}^{-1}(\alpha)\in {\Phi}^-}}|\langle x,\alpha^\vee\rangle-1|+\sum_{\substack{\alpha\in {\Phi}^+\\ {w_J}^{-1}(\alpha)\in {\Phi}^+}}|\langle x,\alpha^\vee\rangle |=l(x)-l(w_J).$$
\end{proof}

{\bf 2.4.} We list some properties of the $a$-function, which were proved in \cite{lusztig1985cells}.
 Let $w, u, v\in W$.

(a) If  $w\underset{LR}\leq u$, then $a(w)\geq a(u)$. In particular, if $w\underset{LR}\sim u$, then $a(w)=a(u)$. In other words,  a-function is constant  on each two-sided cell. Therefore for any integer $i\ge 0$,    $W_{(i)}=\{w\in W\ |\ a(w)=i\}$ is a union of some two-sided cells.

(b) If  $w\underset{L}\leq u$ (respectively,  $w\underset{R}\leq u$; $w\underset{LR}\leq u$) and $a(w)=a(u)$ , then $w\underset{L}\sim u$ (respectively,  $w\underset{R}\sim u$; $w\underset{LR}\sim u$).

(c) For any proper subset $I$ of $S$, if the subgroup $W_I$ generated by $I$ is finite, then for the longest element  $w_I$  of $W_I$,  we have  $a(w_I)=l(w_I)$. If $w=u\cdot w_I\cdot v$, then $a(w)\geq a(w_I)$.

(d) Let ${w_0}$ be the longest element of the Weyl group ${W_0}$,
then for any $w\in W$, $a(w)\leq l({w_0})=a({w_0})=:\nu$, where $\nu$ is the number of positive roots in $\Phi$.

(e) If $w=u\cdot v$, then $w\underset{L}\leq v,\ w\underset{R}\leq u$,
in particular, $a(w)\geq a(v),\ a(u)$.

For $w,u\in W$ we write $w-u$ if $w\le u$ and the Kazhdan-Lusztig polynomial $P_{w,u}$ is of degree $(l(u)-l(w)-1)/2$ or $u\le w$ and the Kazhdan-Lusztig polynomial $P_{u,w}$ is of degree $(l(w)-l(u)-1)/2$.

(f) Let $z, z'\in W$ be such that $z'- z$. If  $\mathscr{R}(z')\not\subset\mathscr{R}(z)$ and $\mathscr{L}(z')\not\subset\mathscr{L}(z)$.  Assume $a(z)<\infty$, then $a(z')>a(z)$.

\medskip

Let $\mathcal{D}=\{w\in {W_a}\ |\ 2 \deg {P_{e,w}}=l(w)-a(w)\}$. Then $\mathcal{D}$ is a finite set and all  elements in $\mathcal D$ are involutions. These involutions are  called distinguished involutions.

Lusztig proved that  $W$ has  only finitely many left cells  and every left cell (respectively, right cell) contains a unique   distinguished involution (see \cite{lusztig1985cells}).

{\bf 2.5.} According to Lemma 4.2 of \cite{nanhua1989approach}, for any $ x\in X,$  there exists  a unique positive integer $n(x)$ such that
\begin{alignat}{2}& x^m\underset{L}{\sim} x^n\quad if \ m,\ n\ge n(x),\\
&x^{n(x)}\underset{L}{\not\sim}  x^{n(x)-1}.\end{alignat}

Since $w\underset{L}\sim u$ if and only if $w^{-1}\underset{R}\sim u^{-1}$, the conclusions above are equivalent to   for any $ x\in X$, there exists  a unique positive integer $n(x)$ such that
\begin{alignat}{2}& x^m\underset{R}\sim x^n\quad if \ m,\ n\ge n(x),\\
&x^{n(x)}\underset{R}{\not\sim}  x^{n(x)-1}.\end{alignat}

It naturally arises a question : $n(x)=?$. There are some evidences showing that $n(x)$ is not big. And for an affine Weyl group of type $\tilde A,\ \tilde B,\ \tilde C,\ \tilde D$, it is possible that $n(x)$ is no larger than 3.

In next section we will show that $n(x)=1$ for any translation in an affine Weyl group of type $\tilde A_{n-1}$. For other types, it is difficult to compute $n(x)$ in general, but  for  translations in the lowest two-sided cell and the second lowest two-sided cell we have similar results.

\section{$n(x)=1$ for  type $\tilde A_{n-1}$ or  the lowest two-sided cell }

\def\ve{\varepsilon}

We keep the notations in section 1. In particular, $W_a$ and $W$ stand for an  affine Weyl group and extended affine Weyl group associated to an irreducible root system $\Phi$ respectively.

{\bf 3.1.} Assume that $W_a$ is of type $\tilde A_{n-1}$. The Coxeter graph of $W_a$ is as below:
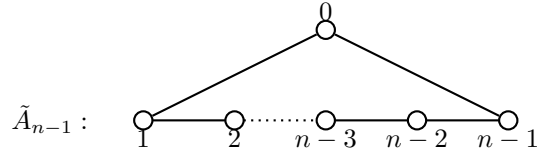
\begin{center}
  \begin{tikzpicture}[scale=.6]
    \draw (-1,0) node[anchor=east]  {$\tilde A_{n-1}:$};
     \draw[thick] (4 cm,2) circle (.2 cm) node [above] {$0$};
    \draw[thick] (0 cm,0) circle (.2 cm) node [below] {$1$};
    \draw[thick] (2 cm,0) circle (.2 cm) node [below] {$2$};
    \draw[thick] (4 cm,0) circle (.2 cm) node [below] {$n-3$};
    \draw[thick] (6 cm,0) circle (.2 cm) node [below] {$n-2$};
    \draw[thick] (8 cm,0) circle (.2 cm) node [below] {$n-1$};
    \draw[thick] (.2 cm,.1) -- (3.8 cm,1.9 cm);
    \draw[thick] (0.2 cm,0) -- +(1.6 cm,0);
    \draw[dotted,thick] (2.2 cm,0) -- +(1.6 cm,0);
    \draw[thick] (4.2 cm,0) -- +(1.6 cm,0);
    \draw[thick] (6.2 cm, 0 cm) -- +(1.6 cm,0);
    \draw[thick] (7.8 cm, 0.1 cm) -- (4.2 cm,1.9 cm);
  \end{tikzpicture}
\end{center}
Correspondingly, the simple refections in $W_a$ are denoted by $s_1,\ ..., s_{n-1},\ s_0$. Let $\Pi=\langle\pi\rangle$ be an infinite cyclic group generated by $\pi$ and let $W=\Pi\ltimes W_a$ such that
 $${s_1}\pi=\pi{s_0},\  {s_2}\pi=\pi{s_1},\  \cdots,\ {s_{0}}\pi=\pi{s_{n-1}}.$$
 Then $W$ is the extended affine Weyl group associated to the general linear group $GL_n(\mathbb C)$. In $W$, the center $Z$ is generated by $\pi^n$. The quotient group $ \tilde W:=W/Z$ is just the extended affine Weyl group associated to a root system of type $\tilde  A_{n-1}$.
There are some advantages to work with this extended affine Weyl group $W$  since it can be realized as a permutation group of $\mathbb Z$. We explain the realization. Let
$${W_*} :=\{\sigma:\mathbb{Z}\to \mathbb{Z}\, |\, \sigma(i+n)=\sigma(i)+n,\  \sum_{i=1}^n(\sigma(i)-i)\equiv 0(\text{mod }n),\ \forall i \in \mathbb{Z}\}.$$
Then we have

(a) (Lusztig)(\cite{xi2002based} Lemma 2.1.1) \label{*} $W$ is isomorphic to $W_*$.

{\bf 3.2.} We shall identify $W$ with $W_*$.  Let $X$ be the subgroup of $W$ consisting of translations, then $X$ is  generated by the following translations $(i=1,\ 2,\ ...,\ n)$:
$$\tau_i(j)=\begin{cases} j+n,&\quad\text{if } j\equiv i(\text{mod}\ n);\\
j,&\quad\text{if } j\not\equiv i(\text{mod}\ n).\end{cases}$$
The simple reflections $s_i$ ($i=0,\ 1,\ 2,\ ..., n-1$) are defined as follows:
$$s_i(j)=\begin{cases} j+1,&\quad\text{if } j\equiv i(\text{mod}\ n);\\
j-1, &\quad\text{if } j\equiv i+1(\text{mod}\ n);\\
j,&\quad\text{if } j\not\equiv i, i+1(\text{mod}\ n).\end{cases}$$
The element $\pi$  is defined by $\pi(i)=i+1$ for all integer $i$.

The simple reflections $s_1,\ ...,s_{n-1}$ generate the Weyl group $W_0$, which is isomorphic to the symmetric group of $n$ letters, and the simple reflections $s_0,\ s_1,\ ...,s_{n-1}$ generate the affine Weyl group $W_a$, which is of type $\tilde A_{n-1}$.

For $1\le i\le n$, set $x_i=\tau_1\tau_2\cdots\tau_i$, then $x_1,\ ...,\ x_{n-1}$ are translations corresponding to fundamental dominant weights $\lambda_1,\ ...,\ \lambda_{n-1}$ respectively. It is easy to see the following relations:

$$\begin{aligned}
&{x_1}=\pi{s_{n-1}}{s_{n-2}}\cdots {s_2}{s_1};\\
&{x_2}=\pi^2({s_{n-2}}{s_{n-3}}\cdots{s_2}{s_1})({s_{n-1}}{s_{n-2}}\cdots{s_2});\\
&\cdots\\
&{x_i}=\pi^i({s_{n-i}}{s_{n-i-1}}\cdots {s_1})({s_{n-i+1}}{s_{n-i}}\cdots{s_2})\cdots({s_{n-1}}{s_{n-2}}\cdots{s_i});\\
&{x_{n-1}}=\pi^{n-1}{s_{1}}{s_{2}}\cdots{s_{n-2}}{s_{n-1}};\\
&{x_{n}}=\pi^{n}.
\end{aligned}$$

One main result of this section is as follows:

{\bf Theorem 3.3.}\label{An} Keep the notations above. In particular,  $W=W_*$ is the extended affine Weyl group associated to $GL_n(\mathbb C)$. For any $x\in X$, $w\in W_0$, we have
\begin{alignat}{2} &x\underset{LR}\sim  wxw^{-1}\\
&x\underset{L}\sim x^m,\quad x\underset{R}\sim x^m,\qquad \text{for any positive integer m}\end{alignat}

{\bf 3.4.} To prove the theorem we recall a description of two-sided cells of $W$ proven by Shi and Lusztig independently.

Let $w\in W$, $j_1,\ j_2,\ ...,\ j_k$ be integers. Call $j_1,\ j_2,\ ...,\ j_k$ a d-antichain of $w$ of length $k$ if the integers satisfy the following relations:

\noindent (1) $j_k-n<j_1<j_2<\cdots<j_k$,

\noindent(2) $w(j_k)-n<w(j_1)<w(j_2)<\cdots w(j_k).$

We call a subset $Y$ of $\mathbb Z$ a d-antichain family set of $w$ of index $q$ if : (1) the elements in $Y$ are noncongruent to each other modulo $n$, (2) $Y$ is a disjoint union of $q$ d-antichains: $A_1,\ A_2,\ ...,\ A_q$. We also call $\{A_1,\ ...,\ A_q\}$ a d-antichain family of $w$.

Let $d_i$ be the maximal one among the cardinalities of all d-antichain  family sets  of $w$ of index $i$. Then $d_1\le d_2\le\cdots\le d_n=n$.  According to Theorem 1.5 in \cite{greene1976some}, we have
$d_1\ge d_2-d_1\ge d_3-d_2\ge\cdots\ge d_n-d_{n-1}$. Then we get a partition of $n$ as
$$\mu(w)=(d_1,\ d_2-d_1,\ \cdots,\ d_n-d_{n-1}).$$

The following result describes the two-sided cells of $W$ completely (\cite{shi1986kazhdan}, \cite{lusztig1985two}).

(a) (Shi,  Lusztig) For any $w,\ u\in W$, then $w\underset{LR}{\sim}u \ \Longleftrightarrow \mu(w)=\mu(u)$.

\medskip

{\bf 3.5.} Now we can prove Theorem 3.3.

(1) Firstly,  We prove for any $w\in W_0$ and $x\in X$,
$$wxw^{-1}\underset{LR}\sim x.$$
Since $W_0$ is generated by $s_1,\ \cdots,\ s_{n-1}$, it is no harm to assume that $w=s_i$ is a simple reflection.
By definition (see subsection 3.2), we have
$$x(j)=j+\theta_jn,\quad \theta_j\in\mathbb Z .$$
Suppose $j_k-n<j_1<j_2<\cdots<j_k$ is a d-antichain of $x$, then
$$j_k+\theta_{j_k}n-n<j_1+\theta_{j_1}n<j_2+\theta_{j_2}<\cdots<j_k+\theta_{j_k}n.$$
Thus $\theta_{j_1}=\theta_{j_2}=\cdots=\theta_{j_k}$.

To consider the antichains of $x$ and $s_ixs_i$, we can only consider the antichains in the set $\{1,\ 2,\ ...,\ n\}$. For any $j$ in this set,
$$s_ixs_i(j)=\begin{cases} x(j)=j+\theta_jn,\quad&\text{if}  \ j\ne i,\ i+1;\\
i+\theta_{i+1}n,&\text{if}  \ j=i;\\
i+1+\theta_{i}n,&\text{if} \ j=i+1.\end{cases}$$

Therefore,

(i) if the d-antichain $ j_1<j_2<\cdots<j_k$ of $x$ in $\{1,\ 2,\ \cdots,\ n\}$ does not contain $i,\ i+1$, then it is also a d-antichain of $s_ixs_i$;

(ii) if the d-antichain $ j_1<j_2<\cdots<j_k$ of $x$ in $\{1,\ 2,\ \cdots,\ n\}$ contains both $i$ and $ i+1$ , then $\theta_i=\theta_{i+1}$, thus it is also a d-antichain of $s_ixs_i$;

(iii) if the d-antichain $ j_1<j_2<\cdots<j_k$ of $x$ in $\{1,\ 2,\ \cdots,\ n\}$ contains $i=j_a $ without $i+1$, we replace $i$ by $i+1$ to get a d-antichain of $s_ixs_i$ whose length is still $k$;

(iv) if the d-antichain $ j_1<j_2<\cdots<j_k$ of $x$ in $\{1,\ 2,\ \cdots,\ n\}$ contains $i+1=j_a  $ without $i$, we replace $i+1$ by $i$ to get a d-antichain of $s_ixs_i$ whose length is still $k$.

In summary, for a d-antichain family set of $x$ of index $q$, we can obtain a d-antichain family set of $s_ixs_i$ of index $q$ correspondingly, and the two sets are of the same cardinality.
As a consequence, we have $\mu(x)$=$\mu(s_ixs_i)$. According to 3.4 (a), $x$ and $s_ixs_i$ are in the same two-sided cell.

(2) Now we prove for any translation $x$ and any positive integer $m$, $x\underset{L} \sim x^m,\quad x\underset{R}\sim  x^m$. In the case of $x$ being a dominant weight, this  result was proved in \cite{xi2002based}. Actually, a stronger result was proved in \cite{xi2002based}:

for $J\subset  \{1, 2, \cdots, n-1\},$ set  $ {x_J}:=\prod\limits_{i\in J}{x_i}$.

(v) if every ${a_i}\geq 1$, then $\prod\limits_{i\in J}{x_i}^{a_i}\underset{R}{\sim} {x_J} $;

(vi) if every ${a_i}\geq 1$, and $k\in\mathbb{Z}$,
then ${x_n}^k\prod\limits_{i\in J}{x_i}^{a_i}\underset{R}{\sim} {x_J} $;

(vii) ${x_1}{x_2}\cdots{x_{n-1}}$ lies in the lowest two-sided cell $W_{(\nu)}$ of $W$.

As is known, there exists $w\in W_0$ such that $wxw^{-1}$ is dominant. Then
$$wxw^{-1}\underset{LR}\sim wx^mw^{-1}.$$
According to the result given in (1), we have $x\underset{LR}{\sim} wxw^{-1}\underset{LR}{\sim} wx^mw^{-1}\underset{LR}{\sim} x^m$. Theorem 3.3 proves valid.\qed

A consequence of Theorem 3.3 is that the  $a$-function value of $x$ equals to the $a$-function value of $x^m$  for any translation $x$ in an extended affine Weyl group of type $\tilde A_{n-1}$.

{\bf Example 3.6.}  Let $W_a$ be an affine Weyl group of type $\tilde G_2$ generated by $s_0,\ s_1,\ s_2$, where $s_1$ is the simple reflection corresponding to the short root. Then $x_1=s_0s_1s_2s_1s_2s_1$ and $x_2=s_0s_1s_2s_1s_0s_2s_1s_2s_1s_2$  and they generate the normal subgroup $X$ of $W$.
It is easy to see
$$x_1\underset{LR}{\not\sim} {x_1}^2,\quad {x_1}^2\underset{L}{\sim} {x_1}^m,\quad x_2\underset{L}{\sim } {x_2}^k, \quad x_1x_2\underset{L}{\sim} {x_1}^j{x_2}^k,$$
where $j,\ k\ge 1$, $m\ge 2$. So for an affine Weyl group of type $\tilde G_2$, there exists $x$ such that $n(x)>1$. Notice that $s_1x_1s_1\underset{LR} {\not\sim }x_1$, thus the assertions of Theorem 3.3 are not true for an affine Weyl group of type $\tilde G_2$.

{\bf 3.7.} In the rest of this section, $W=W_0\ltimes X$ stands for an extended affine Weyl group associated to an irreducible root system $\Phi$. We prove that the assertions of Theorem 3.3 remain valid for the lowest two-sided cell of $W$.

Let $w_0$ be the longest element in $W_0$. Then $\nu=l(w_0)=|\Phi^+|$ is the upper bound of the $a$-function on $W$ (see \cite{lusztig1985cells}).

Recall that $S$ is the set of simple reflections of $W$. For $J\subset S$, the subgroup $W_J$ generated by $J$ is called a parabolic subgroup. If $W_J$ is finite, we denote $w_J$ the longest element in the  subgroup $W_J$.

(a) (Shi, \cite{shi1987two}) The following set is a two-sided cell of $W$, which is the lowest one with respect to the partial order $\underset{LR}\le$ on the set of two-sided cells of $W$:
$$\begin{aligned}
W_{(\nu)}=&\{w\in W\ |\ a(w)=\nu\}\\
=&\{w\cdot w_J\cdot u\ |\ w, \ u\in W, \ W_J\text{ is\ a\ parabolic\ subgroup\ of\ }W,\ l(w_J)=l(w_0)\} .
\end{aligned}$$

It is natural to call $W_{(\nu)}$ the lowest two-sided cell of $W$.
Let $x_I=\prod\limits_{i=1}^{n}x_i$ be the product of all fundamental dominant weights, then $l(x_Is_i)<l(x_I)$ implies that $x_I=w\cdot w_0$ by Lemma 2.3, and we have the following theorem to describe all the translations in $W_{(\nu)}$.

{\bf Theorem 3.8.} Let $x\in X\subset W$. Choose  $w\in W_0$  such that  $wxw^{-1}=x_1^{a_1}x_2^{a_2}\cdots x_n^{a_n}$ is in $X^+$. We have

 (a) $x\in W_{(\nu)}$ if and only if every $a_i$ is a positive integer;

 (b) if $x\in W_{(\nu)}$, then for any $u\in W_0$, we have $uxu^{-1}\in W_{(\nu)}$;

 (c) if $x\in W_{(\nu)}$, then for any positive integer $a$, we have $x^a\in W_{(\nu)}$.

{\bf 3.9.} To prove this theorem, we need the coordinate form for  $w\in W$ which was defined by Jianyi Shi  in \cite{shi1987alcoves}. Be aware that in this article $s_0=t_{-\alpha_0}s_{\alpha_0}$, in \cite{shi1987alcoves} $s_0=s_{\alpha_0}t_{-\alpha_0}$, where $-\alpha_0$ is the highest short root in $\Phi^+$. Since in this article $W$ acts on $E$ on the left and in \cite{shi1987alcoves} $W_a$ (or $W$) acts on $E$ on the right, the alcove corresponding to $w\in W$ in this paper is the same alcove determined by $w^{-1}$ in \cite{shi1987alcoves}. Hence, the coordinate form of $w$ in this article equals to the coordinate form of $w^{-1}$ in \cite{shi1987alcoves}. By definition, the coordinate form of $w$ is a  $\Phi$-tuple   ${(k(w, \alpha))}_{\alpha\in\Phi}$ of integers satisfying the following conditions:

(a) for any $\alpha\in\Phi$, $k(e,\alpha)=0$, where $e$ is the identity element in $W$;

(b) for any $0\leq i\leq n$,
$$k(s_i,\ \alpha)=
\begin{cases}
{0}, &{\alpha\not=\pm{\alpha_i}}\\
{\mp1}, &{\alpha=\pm{\alpha_i}}
\end{cases}$$

(c) if $w'=w\cdot{s_i},\ (0\leq i\leq n)$, then $\forall \alpha\in\Phi$, $$k(w',\ \alpha)=k(w, \ \alpha)+k(s_i, \bar w^{-1} (\alpha)),$$
where $\bar w$ is the image of $w$ under the natural homomorphism $W/X\to W_0$.

For any $\alpha\in\Phi,$ we have $ k(w, \ -\alpha)=-k(w,\  \alpha)$, so the $\Phi$-tuple $(k(w,\ \alpha))_{\alpha\in\Phi}$ is determined by the $\Phi^+$-tuple $(k(w,\  \alpha))_{\alpha\in\Phi^+}$ completely. For convenience we also say that  the $\Phi^+$-tuple $(k(w,\  \alpha))_{\alpha\in\Phi^+}$ is the coordinate form of $w$.

The following result was proved by Jianyi Shi in  \cite{shi1987alcoves} and \cite{shi1987two}.

{\bf Proposition 3.10.} (Shi) Given $w=xu\in W,\ u\in W_0,\ x\in X$, and the corresponding coordinate form  $(k(w,\ \alpha))_{\alpha\in\Phi^+}$, we have:

(a) for any $\alpha\in\Phi^+$, $k(w,\ \alpha)=\langle x, \alpha^\vee\rangle+k(u,\ \alpha)$,
where $$k( u,\ \alpha)=
\begin{cases}
{0}, &{ u^{-1}(\alpha)\in\Phi^+},\\
{-1}, &{ u^{-1}(\alpha)\in\Phi^-}.
\end{cases}$$

(We set $\langle x,\alpha^\vee\rangle$  to be $\langle\lambda,\alpha^\vee\rangle$ if $x=t_\lambda$.)

(b)\ $\mathscr{L}(w)=\{s_j\in S\, |\, k(w, \alpha_j)<0\}$, and
$\mathscr{R}(w)=\{s_j\in S\, |\, k(w,\ {\bar{w}}(\alpha_j))>0\}$, where $\bar w$ is the image of $w$ under the natural homomorphism $W/X\to W_0$.

(c) ${W_{(v)}}:=\{w\in W \,|\, k(w, \ \alpha)\not=0, \ \forall\alpha\in\Phi\}$.

{\bf 3.11.} Now we can prove Theorem 3.8.
By part (a) and (c) of Proposition 3.10 above,  a translation $x\in W$ lies in the lowest two-sided cell ${W_{(v)}}$ if and only if for any $\alpha\in\Phi$, we have $\langle x,\alpha^\vee\rangle \ne 0$. Then Theorem 3.8 follows.

\section{Some translations in the second lowest two-sided cell $\Omega_{qr}$}
{\bf 4.1.} Keep the notations in section 1. In \cite{shi2011second}, Jianyi Shi studied the second lowest two-sided cell of $W_a$, which has little difference with the second lowest two-sided cell of $W$. Following Shi, we denote $\Omega_{qr}$ the second lowest two-sided cell of $W$. Shi conjectured that $\Omega_{qr}$ contains $|W_0|/2$ left cells and proved that the number of left cells in $\Omega_{qr}$ is at most $|W_0|/2$.

In \cite{shi2011second}, Jianyi Shi gave the table of $a$-function values of $\Omega_{qr}$ as below:

\textbf{Table 1}

\begin{tabular}{cccccccccc}
\hline
      $W_a$ & $\tilde{A_n}$ & $\tilde{C_n}$ &$\tilde{B_n}$ & $\tilde{D_n}$ & $\tilde{E_6}$ & $\tilde{E_7}$ & $\tilde{E_8}$ & $\tilde{F_4}$ & $\tilde{G_2}$\\
\hline
$a(\Omega_{qr})$ & $\frac{1}{2}(n^2-n)$ & $n^2-2n+2$ & $n^2-n$ & $n^2-3n+3$ & 25 &46 & 91 & 16 &3\\
\hline
\end{tabular}

\bigskip

The main aim of this section is to determine some translations in the second lowest two-sided cell $\Omega_{qr}$. We hope to use the translations to understand Shi's conjecture on the number of left cells in $\Omega_{qr}$.

Let $x\in X$ be a translation in $W$. For convenience, we set $\langle x,\alpha^\vee\rangle=\langle \lambda,\alpha^\vee\rangle$ if $x=t_\lambda$ and $\alpha\in\Phi$.  By   Proposition 3.10 (c),  if $x\in X\cap\Omega_{qr}$, then there  exists some $\alpha\in\Phi$ such that $\langle x,\alpha^\vee\rangle=0$. For $J\subset I=\{1,2,...,n\}$, set $x_J=\prod_{i\in J}x_i$. Theorem 4.1 in \cite{nanhua1989approach} says that there exists some positive integer $r$ and $1\le i\le n$ satisfying
$${x_{I-\{i\}}}^r=\prod_{j\ne i, j\in I}{x_j}^r\in \Omega_{qr}.$$

To go further we firstly establish a technical result.

{\bf Lemma 4.2.} Assume that the root system $\Phi$ has rank $n$. Let $s_1,\ \cdots,\ s_n$ be the simple reflections in $W_0$ and  $x_1,\ \cdots, \ x_n$ are the corresponding translations in $W$ (see subsection 1.3). Set $I_k=\{1,\ \cdots,\ n\}\backslash\{k\}$, $1\le k\le n$. Then
$$x_{I_k}s_k=u\cdot w_0$$
lies in the lowest two-sided cell ${W_{(v)}}$ of $W$,  where $\displaystyle x_{I_k}=\prod_{i\in I_k}x_i$, $w_0$ is the longest element in the  Weyl group $W_0$.

\begin{proof}
By the length formula (1) of Iwahori-Matsumoto,
$$l(x_{I_k}s_k)=|\langle x_{I_k}, {\alpha_k}^\vee\rangle-1 |+\sum_{\substack{\alpha\in{\Phi}^+\\ \alpha\not=\alpha_k}}|\langle x_{I_k},\alpha^\vee\rangle|=l(x_{I_k})+1,$$
$$l(x_{I_k}s_kw_0)=\sum_{\substack{\alpha\in {\Phi}^+\\ \alpha\not=\alpha_k}}|\langle x_{I_k},\alpha^\vee\rangle-1|+|\langle x_{I_k}, {\alpha_k}^\vee\rangle |=l(x_{I_k})+1-l(w_0),$$
then, $l(x_{I_k}s_k)=l(x_{I_k}s_kw_0)+l(w_0)$. Consequently, $x_{I_k}s_k=x_{I_k}s_kw_0\cdot w_0$.
\end{proof}

{\bf Proposition 4.3.} Assume that  $x\in X\cap  \Omega_{qr}$, then
for any positive integer $k$, we have $x^k\in\Omega_{qr}$ and $x^k\underset{L}\sim x, \ x^k\underset{R}\sim x.$

\begin{proof}
According to  Proposition 3.10 (c),  there exists some $\alpha\in\Phi$ such that $\langle x,\alpha^\vee\rangle=0$. This implies that  $\langle  x^k, \alpha\rangle=k\cdot 0=0$ for any positive integer $k$. Hence $a(x^k)\leq a(\Omega_{qr})$.

On the other hand, by 2.4(e) we have $a(x^k)\geq a(x)=a(\Omega_{qr})$. This forces that $a(x^k)= a(x)$. Using 2.4(e) we see that $x^a\underset{L}\leq x$ and $x^k \underset{R}\leq  x$. Applying 2.4(b) we get $x^k\underset{L}\sim x,\  x^k\underset{R}\sim x.$
\end{proof}

Then we have $n(x)=1$ for any translation $x\in\Omega_{qr}$.
\medskip

{\bf 4.4.} Now we are going to give a family of translations in the second lowest two-sided cell $\Omega_{qr}$ of the extended affine Weyl group $W$ case by case. Note that the root system $\Phi$ has rank $n$. The simple reflections $s_1,\ ...,\ s_n$  in $W_0$ are numbered as usual. The corresponding fundamental dominant weights in $X$ are denoted by $x_1,\ ..., x_n$. The simple reflection in $W\setminus W_0$ is denoted by $s_0$.

 For $1\le k\le n$ we set $I_k=\{1,\  2,\ ...,\ n\}\setminus\{k\}$ and $x_{I_k}=\prod_{i\in I_k}x_i$.

 {\bf Theorem 4.5.} Assume that $W$ is of type {\bf $\tilde A_n$ or $ \tilde D_n \ (n\ge 4)$ or $\tilde E_n\ (n=6,7,8)$}. Then a translation $x\in X^+$ is in $\Omega_{qr}$  if and only if there exists $1\le k\le n$ and positive integers $a_i$ for all $i\in I_k$ such that
  $$x=\prod_{i\in I_k} x_i^{a_i}.$$

\begin{proof}
Suppose that $x=\prod_{i\in I_k} x_i^{a_i}$,  $ a_i\ge 1$ for all $i\in I_k$. Then  $\langle x, {\alpha_k}^\vee\rangle=0$. By the proof of Theorem 3.8, we know that $x$ is not in the lowest two-sided cell $W_{(\nu)}$ of $W$. It is sufficient to show $\displaystyle x_{I_k}= \prod_{i\in I_k} x_i\in\Omega_{qr}$ since $x\underset{LR}\sim x_{I_k}$.

By Lemma 4.2, $x_{I_k}s_k=u\cdot w_0$, where $w_0$ is the longest element of the  Weyl group $W_0\subset W$. Thus $x_{I_k}=u\cdot (w_0s_k)$. According to the proof of Proposition 5.7 in \cite{shi2011second} , for an affine Weyl group of type $\tilde A_n,\ \tilde D_n \ (n\ge 4),\ \tilde E_n\ (n=6,7,8)$, $w_0s_k\in\Omega_{qr}$. This forces $x_{I_k}\in\Omega_{qr}$.

 Now suppose that   $x=\prod_{i=1 }^n x_i^{a_i}\in X^+$ is in $\Omega_{qr}$. By Theorem 3.8 (a), there exists some $1\le k\le n$ such that $a_k=0$.  If in addition,  $a_j=0$ for some $j\ne k$, then $s_jx=xs_j$ and this implies that $a(x)< a(s_jx)$ by 2.4(f). Since  $k(s_jx, \alpha_k)=0$, by Proposition 3.10, we see   $s_jx\not\in W_{(\nu)}$. Thus $a(x)<a(s_jx)\leq a(\Omega_{qr})$ and this contradicts $x\in \Omega_{qr}$.  Therefore there is only one $k$ such that $a_k=0$ and all other $a_i$ are positive integers.
 \end{proof}

{\bf Corollary 4.6.}  Assume that $W$ is of type $\tilde A_n$. Then $x\in X$ is in $\Omega_{qr}$ if and only if there exists some $w\in W_0$ and $1\le k\le n$ such that $wxw^{-1}=\prod_{i\in I_k} x_i^{a_i}, \ a_i\ge 1$.

\begin{proof}
It follows from Theorem 4.5 and Theorem 3.3 directly.
\end{proof}

\def\ve{\varepsilon}

\medskip

{\bf 4.7.} Now we deal with type  $\tilde B_n$ $(n\ge 3)$.

 Let $E$ be the n-dimensional  euclidean space with the standard orthogonal basis $\ve_1,\ \ve_2,\ ...,\ \ve_n$, $\Phi=\{\pm(\ve_i\pm\ve_j),\ \pm2\ve_k\,|\, 1\le i<j\le n,\ 1\le k\le n\}.$  The extended affine Weyl group $W$ associated to $\Phi$ is of type $\tilde B_n$. The Coxeter graph of $W$ is as follows ( see \cite{shi2011second} ):
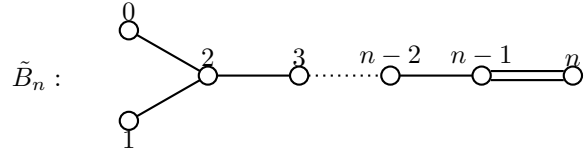
\begin{center}
  \begin{tikzpicture}[scale=.6]
    \draw (-1,0) node[anchor=east]  {$\tilde B_n:$};
    \draw[thick] (2 cm,0) circle (.2 cm) node [above] {$2$};
    \draw[xshift=2 cm,thick] (150:2) circle (.2 cm) node [above] {$0$};
    \draw[xshift=2 cm,thick] (210:2) circle (.2 cm) node [below] {$1$};
    \draw[thick] (4 cm,0) circle (.2 cm) node [above] {$3$};
    \draw[thick] (6 cm,0) circle (.2 cm) node [above] {$n-2$};
    \draw[thick] (8 cm,0) circle (.2 cm) node [above] {$n-1$};
    \draw[thick] (10 cm,0) circle (.2 cm) node [above] {$n$};
    \draw[xshift=2 cm,thick] (150:0.2) -- (150:1.8);
    \draw[xshift=2 cm,thick] (210:0.2) -- (210:1.8);
    \draw[thick] (2.2,0) --+ (1.6,0);
    \draw[dotted,thick] (4.2,0) --+ (1.6,0);
    \draw[thick] (6.2,0) --+ (1.6,0);
    \draw[thick] (8.2,0.1) --+ (1.6,0);
    \draw[thick] (8.2,-0.1) --+ (1.6,0);
  \end{tikzpicture}
\end{center}
 \def\ve{\varepsilon}

Fix a simple root system as
$$\Delta=\{\ve_1-\ve_2,\ \ve_2-\ve_3,\ ...,\ \ve_{n-1}-\ve_n,\ 2\ve_n\}.$$
Then the set of positive roots is
$$  \Phi^+=\{ \ve_i\pm\ve_j,\ 2\ve_k\,|\, 1\le i<j\le n,\ 1\le k\le n\},$$
and the set of negative roots is $\Phi^-=-\Phi^+.$

Let\begin{align*}
 \alpha_1&=\ve_1-\ve_2,\ \alpha_2=\ve_2-\ve_3,\ ...,\ \alpha_{n-1}=\ve_{n-1}-\ve_n,\ \alpha_n=2\ve_n;\\
\alpha_{ij}&=\ve_i-\ve_j,\quad \beta_{ij}=\ve_i+\ve_j,\quad1\le i<j\le n,\\
\gamma_k&=2\ve_k,\qquad 1\le k\le n.
\end{align*}
 Then
 \begin{align*}
\alpha_{ij}&=\alpha_i+\alpha_{i+1}+\cdots+\alpha_{j-1}, \quad 1\le i<j\le n;\\
\beta_{ij}&=\alpha_i+\cdots+\alpha_{j-1}+2\alpha_j+\cdots+2\alpha_{n-1}+\alpha_n,\quad 1\le i<j\le n-1;\\
\beta_{in}&=\alpha_i+\alpha_{i+1}+\cdots+\alpha_n,\quad 1\le i\le n-1;\\
\gamma_{k}&=2\alpha_k+\alpha_{k+1}+\cdots+2\alpha_{n-1}+\alpha_n,\quad 1\le k\le n-1;\\
\gamma_n&=\alpha_n.
\end{align*}

   The Weyl group of the root system $\Phi$ is generated by the simple reflections $s_1,\ ...,\  s_n$. The corresponding fundamental dominant weights in $X\subset W$ are denoted by $x_1,\ ...,\ x_n$. As usual, $s_0$ is the simple reflection out of $W_0$.

{\bf Theorem 4.8.} Keep the notations in subsection 4.7, so $W$ is an extended affine Weyl group of type $\tilde B_n$. Let
$$x=\displaystyle\prod_{i=1}^{n-1} x_i^{a_i}, \qquad a_1,\ ...,\ a_{n-2}\ge 1,\ a_{n-1}\ge 2.$$
 Then $x\in \Omega_{qr}$.

{\bf 4.9.} Now we give a proof of Theorem 4.8. By Theorem 3.8 (a), $x$ is not in the lowest two-sided cell $W_{(\nu)}$ of $W$. Thus $a(x)\le n^2-n=a(\Omega_{qr})$. Let $J=\{s_0,\ s_1,\ ...,\ s_{n-1}\}$, then $W_J=\langle s_i\,|\, i\in J\rangle$ is  a Weyl group of type $D_n$ and its  longest element $w_J$ is of length $n^2-n$. Therefore to prove the theorem it suffices to show
$$x_1\cdots x_{n-2}x_{n-1}^2=w\cdot w_J\cdot u.$$

According to Lemma  2.3, $\mathscr{R}(x_1\cdots x_{n-2}x_{n-1})=\{s_1,\ s_2,\ ...,\ s_{n-1}\}$, thus
$$x_1\cdots x_{n-2}x_{n-1}=v\cdot w_K,$$
where $K=\{1,\ ...,\ n-1\}$ and $w_K$ is the longest element of $W_K=\langle s_i\,|\, i\in K\rangle$

It suffices to prove
$$w_Kx_{n-1}=w_J\cdot u.$$
We know that $\{s_1,\ ...,\ s_{n-1}\}\subset \mathscr{L}(w_Kx_{n-1})$. It will be done if  we can show $s_0w_Kx_{n-1}\le w_Kx_{n-1}$, or equivalently, $x_{n-1}^{-1}w_Ks_0\le x_{n-1}^{-1}w_K$.

We need the following two formulas:

(a) $l(x_{n-1}^{-1})=n^2-1$,

(b)  $l(x^{-1}_{n-1}w_K)=n^2-1+\frac{n(n-1)}2$.

We  use the length formula (1) of Iwahori-Matsumoto to prove (a) and (b).

Recall the weight lattice  of $\Phi$ is  $P=\{\lambda\in E\,|\, \langle \lambda,\ \alpha^\vee\rangle  \in\mathbb Z,\text{ for\  any\ } \alpha\in \Phi\},$
where $\langle \cdot, \cdot\rangle $ is the inner product in $E$,
$\alpha^\vee=2\alpha/\langle \alpha,\alpha).$
The fundamental dominant weights   $\lambda_1,\ ...,\ \lambda_n$ are defined by  $\langle \lambda_i,\ \alpha_j^\vee\rangle  =\delta_{ij},\  1\le i,j\le n.$

\def\ve{\varepsilon}
By the length formula (1) of Iwahori-Matsumoto, we have
$$l(x_{n-1}^{-1})=\sum_{\alpha\in\Phi^+}|\langle-\lambda_{n-1},\ \alpha^\vee\rangle  |=
\sum_{\alpha\in\Phi^+}|\langle \lambda_{n-1},\ \alpha^\vee\rangle  |.$$
Note that $\alpha_{ij}^\vee=\alpha_{ij},\ \ \beta_{ij}^\vee=\beta_{ij},\ \ \gamma_k^\vee=\frac12\gamma_k$. Hence
\begin{align*}
\langle \lambda_{n-1},\ \alpha^\vee_{ij}\rangle&=\langle \lambda_{n-1},\ \alpha _{ij}\rangle=
\begin{cases}
{0},\quad &{1\le i< j\le n-1}\\
{1},\quad &{1\le i<j=n}
\end{cases};\\
\langle \lambda_{n-1},\ \beta^\vee_{ij}\rangle&=\langle \lambda_{n-1},\ \beta _{ij}\rangle=
\begin{cases}
{2},\quad &{1\le j\le n-1}\\
{1},\quad &{1\le i<j=n}
\end{cases};\\
\langle \lambda_{n-1},\ \gamma^\vee_{k}\rangle&=\langle \lambda_{n-1},\ \frac{\gamma _{k}}{2}\rangle=
\begin{cases}
{1},\quad &{1\le k\le n-1}\\
{0},\quad &{k=n}
\end{cases}.
\end{align*}
Thus,
\begin{align*} l(x_{n-1}^{-1})&=\sum_{1\le i<j\le n}|\langle \lambda_{n-1},\ \alpha^\vee_{ij}\rangle|+\sum_{1\le i<j\le n}|\langle \lambda_{n-1},\ \beta^\vee_{ij}\rangle|+\sum_{1\le k\le n}|\langle \lambda_{n-1},\ \gamma^\vee_{k}\rangle|\\
&=n-1+(2\cdot\frac{(n-1)(n-2)}2+n-1)+n-1=n^2-1.
\end{align*}
Formula (a)  is proved. Formula (b) follows from $l(x_{n-1}^{-1}w_K)=l(x_{n-1}^{-1})+l(w_K)$.

\medskip
\def\ve{\varepsilon}
Now we are able to show $x_{n-1}^{-1}w_Is_0\le x_{n-1}^{-1}w_I$, which is equivalent to show $l(x_{n-1}^{-1}w_Is_0)=n^2-2+\frac{n(n-1)}2$.

By definition, $s_0=t_{-\alpha_0}s_{\alpha_0}$, where $-\alpha_0=\ve_1+\ve_2=\lambda_2$ is the highest short root. Thus
$$x_{n-1}^{-1}w_Ks_0=x_{n-1}^{-1}w_Kt_{\lambda_2}s_{\alpha_0}=x_{n-1}^{-1}t_{w_K(\lambda_2)}w_Ks_{\alpha_0}.$$

We compute $w_K(\lambda_2)$. Notice that $s_i(\lambda_2)=\lambda_2$ if $i\ne 2$.  Set $L=\{1,\ 3,\ \cdots,\ n-1\}$, $w_L$ is the longest element of $W_L=\langle s_k\,|\, k\in L\rangle$. We have
$$w_K=s_{n-2}s_{n-3}\cdots s_2s_1s_{n-1}s_{n-2}\cdots s_3s_2w_L,\quad w_L(\lambda_2)=\lambda_2.$$
 Then
\begin{align*}
 w_K(\lambda_2)&=s_{n-2}s_{n-3}\cdots s_2s_1s_{n-1}s_{n-2}\cdots s_3s_2(\lambda_2)\\
&=s_{n-2}s_{n-3}\cdots s_2s_1(\lambda_2-\alpha_2-\cdots-\alpha_{n-1})\\
&=s_{n-2}s_{n-3}\cdots s_2s_1(\lambda_2-\alpha_1-\alpha_2-\cdots-\alpha_{n-1})\\
&= \lambda_2-\alpha_1-2\alpha_2-\cdots-2\alpha_{n-2}-\alpha_{n-1}\\
&=\alpha_1+2\alpha_2+\cdots+2\alpha_{n-1}+\alpha_n-\alpha_1-2\alpha_2-\cdots-2\alpha_{n-2}-\alpha_{n-1}\\
&=\alpha_{n-1}+\alpha_n=\ve_{n-1}-\ve_n+2\ve_n=\ve_{n-1}+\ve_n.
\end{align*}

Therefore, by the length formula (1) of Iwahori-Matsumoto,  we have
\begin{align*}
 &\ \qquad  l(x_{n-1}^{-1}w_Is_0)=l(x_{n-1}^{-1}t_{\ve_{n-1}+\ve_n}w_Is_{\alpha_0})\\
&=\sum_{\substack{\alpha\in{\Phi^+}\\ s_{\alpha_0}w_I(\alpha)\in{\Phi^+}}}|\langle -\lambda_{n-1}+\ve_{n-1}+\ve_n,\ \alpha^\vee\rangle  |+
\sum_{\substack{\alpha\in{\Phi^+}\\ s_{\alpha_0}w_I(\alpha)\in{\Phi^-}}}|\langle -\lambda_{n-1}+\ve_{n-1}+\ve_n,\ \alpha^\vee\rangle  -1|.
\end{align*}

We need to figure out the positive roots mapped by $s_{\alpha_0}w_I$ into $\Phi^+$ or $\Phi^-$ respectively in detail.

We write $\alpha>0$ if  $\alpha$ is a positive root and write $\alpha<0$ if $\alpha$ is a negative root.
Assume $1\le i<j\le n,\ 1\le k\le n$, it is easy to get
\begin{align*}
 s_{\alpha_0}w_I(\alpha_{ij})&=s_{\alpha_0}w_I(\ve_i-\ve_j)=s_{\alpha_0}(\ve_j-\ve_i)\\
&=\begin{cases}
{ \ve_j-\ve_i<0},\quad &{i,\  j\not\in\{1,\ 2\}}\\
{\ve_2-\ve_1<0},\quad &{(i,\ j)=(1,\ 2)}\\
{\ve_2+\ve_j>0},\quad &{i=1,\ j>2}\\
{\ve_1+\ve_j>0},\quad &{i=2,\ j>2}.
\end{cases}
\end{align*}

\begin{align*}
s_{\alpha_0}w_I(\beta_{ij})&=s_{\alpha_0}w_I(\ve_i+\ve_j)=s_{\alpha_0}(\ve_{n-i+1}+\ve_{n-j+1})\\
&=
\begin{cases}
{ \ve_{n-i+1}+\ve_{n-j+1}>0},\quad &{i,\ j\not\in\{n-1,\ n\}}\\
{-\ve_{ 1}-\ve_{2}<0},\quad &{(i,\ j)=(n-1,\ n)}\\
{-\ve_2+\ve_{n-i+1}<0},\quad &{1\le i\le n-2,\ j=n}\\
{-\ve_1+\ve_{n-i+1}<0},\quad &{1\le i\le n-2,\ j=n-1}.
\end{cases}
\end{align*}
\begin{align*}
s_{\alpha_0}w_I(\gamma_k)&=s_{\alpha_0}w_I(2\ve_k)=s_{\alpha_0}(2\ve_{n-k+1} )\\
&=\begin{cases}
{ 2\ve_{n-k+1} >0},\quad &{1\le k\le n-2}\\
{-2\ve_{ 2} <0},\quad &{k=n}\\
{- 2\ve_1 <0},\quad &{k=n-1}.
\end{cases}
\end{align*}

Now we compute $|\langle -\lambda_{n-1}+\ve_{n-1}+\ve_n,\ \alpha^\vee\rangle  |$ and
$|\langle -\lambda_{n-1}+\ve_{n-1}+\ve_n,\ \alpha^\vee\rangle  -1|$ in accordance with
$s_{\alpha_0}w_I(\alpha)>0$ or $<0$.

Case 1: $\alpha\in\Phi^+,\ s_{\alpha_0}w_I(\alpha)>0,.$

(1) When $\alpha=\alpha_{1j}=\ve_1-\ve_j,\ \ j\ge 3$, we have
$$|\langle -\lambda_{n-1}+\ve_{n-1}+\ve_n,\ \alpha^\vee\rangle  |=
\begin{cases}
{ 0},\quad &{1\le j\le n-2}\\
{1},\quad &{j=n-1}\\
{2},\quad &{j=n}.
\end{cases}$$
Thus, these positive roots contribute 3 to the sum in the length formula.

(2) When $\alpha=\alpha_{2j}=\ve_2-\ve_j,\ \ j\ge 3$, we have
$$|\langle -\lambda_{n-1}+\ve_{n-1}+\ve_n,\ \alpha^\vee\rangle  |=
\begin{cases}
{ 0},\quad &{1\le j\le n-2}\\
{1},\quad &{j=n-1}\\
{2},\quad &{j=n}.
\end{cases}$$
Thus,  these positive roots contribute 3 to the sum in the length formula.

(3) When  $\alpha=\beta_{ij}=\ve_i+\ve_j,\ \ i,j\not\in\{n-1,\ n\}$, we have
$$|\langle -\lambda_{n-1}+\ve_{n-1}+\ve_n,\ \alpha^\vee\rangle  |=2 .$$
Thus,  these positive roots contribute $2\cdot \frac{(n-2)(n-3)}2=(n-2)(n-3)$ to the sum in the length formula.

(4) When  $\alpha=\gamma_{k}=2\ve_k,\ \ k\ne n-1,\ n$, we have
$$|\langle -\lambda_{n-1}+\ve_{n-1}+\ve_n,\ \alpha^\vee\rangle  |=1. $$
Thus,  these positive roots contribute $n-2$ to the sum in the length formula.

Case 2: $\alpha\in\Phi^+,\ s_{\alpha_0}w_I(\alpha)<0$,

(1) When $\alpha=\alpha_{ij}=\ve_1-\ve_j,\ \ i,\ j\not\in\{1,2\}$, we have
$$|\langle -\lambda_{n-1}+\ve_{n-1}+\ve_n,\ \alpha^\vee\rangle  -1|=
\begin{cases}
{1},\quad &{3\le i<j\le n-2}\\
{3},\quad &{3\le i\le n-2,\ j=n}\\
{2},\quad &{i=n-1,\ j=n}\\
{2},\quad &{3\le i\le n-2,\ j=n-1}.
\end{cases}$$
Thus, these positive roots contribute $\frac{(n-4)(n-5)}2+3(n-4)+2+2(n-4)=5n-18+\frac{(n-4)(n-5)}2$ to the sum in the length formula.

(2) When $\alpha=\ve_1-\ve_2$, we have
$$|\langle -\lambda_{n-1}+\ve_{n-1}+\ve_n,\ \alpha^\vee\rangle  -1|=1.$$
Thus, this  positive root contributes 1 to the sum in the length formula.

(3) When $\alpha=\ve_{n-1}+\ve_{n}$, we have
$$|\langle -\lambda_{n-1}+\ve_{n-1}+\ve_n,\ \alpha^\vee\rangle  -1|=0.$$
Thus, this  positive root contributes 0 to the sum in the length formula.

(4)  When $\alpha=\beta_{in}=\ve_i+\ve_n,\ \ 1\le i\le n-2$, we have
$$|\langle -\lambda_{n-1}+\ve_{n-1}+\ve_n,\ \alpha^\vee\rangle  -1|=1.$$
Thus, these positive roots contribute $n-2$ to the sum in the length formula.

(5) When $\alpha=\beta_{i,n-1}=\ve_i+\ve_{n-1},\ \ 1\le i\le n-2$, we have
$$|\langle -\lambda_{n-1}+\ve_{n-1}+\ve_n,\ \alpha^\vee\rangle  -1|=2.$$
Thus, these positive roots contribute $2(n-2)$ to the sum in the length formula.

(6) When $\alpha=\gamma_n=2\ve_n$, we have
$$|\langle -\lambda_{n-1}+\ve_{n-1}+\ve_n,\ \alpha^\vee\rangle  -1|=0 .$$
Thus, this  positive root contributes 0 to the sum in the length formula.

(7) When $\alpha=\gamma_{n-1}=2\ve_{n-1}$, we have
$$|\langle -\lambda_{n-1}+\ve_{n-1}+\ve_n,\ \alpha^\vee\rangle  -1|=1 .$$
Thus, this  positive root contributes 1 to the sum in the length formula.

Add up all the numbers in Case 1 and Case 2, finally we get:
$$\begin{aligned}
l(x_{n-1}^{-1}w_Is_0)=&3+3+(n-2)(n-3)+(n-2)\\
&\ +5n-18+\frac{(n-4)(n-5)}2+1+0+(n-2)+2(n-2)+0+1\\
=&n^2-2+\frac{n(n-1)}2=l(x_{n-1}^{-1}w_I)-1.
\end{aligned}$$
which is exactly what we need, and the theorem follows.\qed

\begin{remark}
It is difficult to prove the theorem by using the reduced expressions of the fundamental dominant weights directly. Actually Jinhao Guo discussed the reduced expressions in   section 3.2.2 in his PhD thesis\cite{guo2019kuozhang}, the expressions seem complicated for application except for a few cases. The reduced expressions obtained by Guo are deduced from the formulas in \cite{lusztig1983some}:
$$\begin{aligned}
 &x_1=\tau s_1 s_2s_3\cdots s_{n-1}s_ns_{n-1}\cdots s_3s_2s_1,\\
&x_2=s_1x_1s_1x_1,\\
&x_3=s_2s_1x_1s_1s_2x_2,\\
&.........,\\
&x_{n-1}=s_{n-2}\cdots s_2s_1x_1s_1s_2\cdots s_{n-2}x_{n-2},\\
&x_n=s_{n-1}s_{n-2}\cdots s_2s_1x_1s_1s_2\cdots s_{n-2}s_{n-1}x_{n-1},
\end{aligned}$$

where $\tau\in \Omega$ satisfies $\tau s_0=s_1\tau,\ \tau s_i=s_i\tau$, for $ i=2,\ 3,\ \cdots,\ n$.
\end{remark}

\def\ve{\varepsilon}
\bigskip

{\bf 4.10.} Now we deal with type  $\tilde C_n$ $(n\ge 2)$.

 Let $E$ be the n-dimensional  euclidean space with the standard orthogonal basis $\ve_1,\ \ve_2,\ ...,\ \ve_n$, $\Phi=\{\pm(\ve_i\pm\ve_j),\  \ve_k\,|\, 1\le i<j\le n,\ 1\le k\le n\}.$ The extended affine Weyl group associated to $\Phi$ is of type $\tilde C_n\ (n\ge 2)$.
The Coxeter graph of $W$ is as follows ( see \cite{shi2011second}):
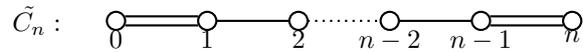
\begin{center}
  \begin{tikzpicture}[scale=.6]
    \draw (-1,0) node[anchor=east]  {$\tilde{C_n}:$};
    \draw[thick] (0 cm,0) circle (.2 cm) node [below] {$0$};
    \draw[thick] (2 cm,0) circle (.2 cm) node [below] {$1$};
    \draw[thick] (4 cm,0) circle (.2 cm) node [below] {$2$};
    \draw[thick] (6 cm,0) circle (.2 cm) node [below] {$n-2$};
    \draw[thick] (8 cm,0) circle (.2 cm) node [below] {$n-1$};
    \draw[thick] (10 cm,0) circle (.2 cm) node [below] {$n$};
    \draw[thick] (.2 cm,.1) -- +(1.6 cm,0);
    \draw[thick] (.2 cm,-.1) -- +(1.6 cm,0);
    \draw[thick] (2.2 cm,0) -- +(1.6 cm,0);
    \draw[dotted,thick] (4.2 cm,0) -- +(1.6 cm,0);
    \draw[thick] (6.2 cm,0) -- +(1.6 cm,0);
    \draw[thick] (8.2 cm, .1 cm) -- +(1.6 cm,0);
    \draw[thick] (8.2 cm, -.1 cm) -- +(1.6 cm,0);
  \end{tikzpicture}
\end{center}

Fix a simple root system as
$$\Delta=\{\ve_1-\ve_2,\ \ve_2-\ve_3,\ ...,\ \ve_{n-1}-\ve_n,\  \ve_n\}.$$
Then the set of positive roots is
$$  \Phi^+=\{ \ve_i\pm\ve_j,\  \ve_k\,|\, 1\le i<j\le n,\ 1\le k\le n\},$$
and the set of negative roots is $\Phi^-=-\Phi^+.$ Let
$$\begin{aligned}
 \alpha_1&=\ve_1-\ve_2,\ \alpha_2=\ve_2-\ve_3,\ ...,\ \alpha_{n-1}=\ve_{n-1}-\ve_n,\ \alpha_n= \ve_n;\\
\alpha_{ij}&=\ve_i-\ve_j,\quad \beta_{ij}=\ve_i+\ve_j,\quad1\le i<j\le n,\\
\gamma_k&=\ve_k,\qquad 1\le k\le n.
\end{aligned}$$
Then
$$\begin{aligned}
 \alpha_{ij}&=\alpha_i+\alpha_{i+1}+\cdots+\alpha_{j-1}, \quad 1\le i<j\le n;\\
\beta_{ij}&=\alpha_i+\cdots+\alpha_{j-1}+2\alpha_j+\cdots+2\alpha_{n-1}+2\alpha_n,\quad 1\le i<j\le n;\\
\gamma_{k}&= \alpha_k+\alpha_{k+1}+\cdots+ \alpha_{n-1}+\alpha_n,\quad 1\le k\le n.
\end{aligned}$$
The fundamental dominant weights $\lambda_1,\ \cdots,\ \lambda_n$ in the weight lattice $P$ of $\Phi$ are well known:
$$\begin{aligned}
 \lambda_i&=\alpha_1+2\alpha_2+\cdots+(i-1)\alpha_{i-1}+i(\alpha_i+\alpha_{i+1}+\cdots+\alpha_n),
\quad 1\le i\le n-1,\\
\lambda_n&=\frac12(\alpha_1+2\alpha_2+\cdots+n\alpha_n).
\end{aligned}$$
In terms of the standard orthogonal basis we have
$$\begin{aligned}
\lambda_i&=\ve_1+\ve_2+\cdots+\ve_i, \quad 1\le i\le n-1,\\
\lambda_n&=\frac12(\ve_1+\ve_2+\cdots+\ve_n).
\end{aligned}$$

Note that $t_{\lambda_i}=x_i$. Using the length formula (1) of Iwahori-Matsumoto, we get

 (a) $l(x_i)=i(2n-i+1),\quad 1\le i\le n-1$;\

 (b)  $l(x_n)=\frac{n(n+1)}2$.

\medskip

The following identities can be found in \cite{lusztig1983some}:
$$\begin{aligned} 
&x_1=s_0 s_1 s_2s_3\cdots s_{n-1}s_n s_{n-1}\cdots s_3s_2s_1\\
&x_2=s_1x_1s_1x_1\\
&x_3=s_2s_1x_1s_1s_2x_2\\
&.........\\
&x_{n-1}=s_{n-2}\cdots s_2s_1x_1s_1s_2\cdots s_{n-2}x_{n-2}\\
&x_n=\tau (s_n s_{n-1} \cdots s_2s_1)(s_{n }s_{n-1}\cdots s_2)\cdots(s_ns_{n-1})(s_n),\end{aligned}$$
where $\tau\in \Omega$ satisfies $\tau s_0=s_n\tau,\ \tau s_i=s_{n-i}\tau$, for $ i=2,3,...,n-1$.

We can determine reduced expressions of the fundamental dominant weights by these formulas.  Jinhao Guo has obtained a reduced expression for each $x_i$ (see Lemma 3.3.2 in \cite{guo2019kuozhang}). But those reduced expressions are not easy to use for the discussion here. Hence, we display a different reduced expression for these $x_i$.

\medskip

{\bf Proposition 4.11.} Keep the notations in subsection 4.10. Then $x_1,\ ...,\ x_n$  are fundamental dominant weights in the extended affine Weyl group $W$ of type $\tilde C_n$. We claim that the following  expressions are reduced:
$$\begin{aligned}
&x_1=s_0s_1s_2\cdots s_{n-1}s_ns_{n-1}\cdots s_2s_1;\\
&x_2=s_0s_1s_2\cdots s_{n-1}s_ns_{n-1}\cdots s_2 s_0s_1s_2\cdots s_{n-1}s_ns_{n-1}\cdots s_2 =(s_0s_1s_2\cdots s_{n-1}s_ns_{n-1}\cdots s_2)^2; \\
&............\\
&x_i=(s_0s_1s_2\cdots s_{n-1}s_ns_{n-1}\cdots s_{i})^i,\quad 1\le i\le n-1;\\
&x_n=\tau (s_{n }s_{n-1}\cdots s_2s_1)(s_ns_{n-1}\cdots s_2)\cdot\cdots\cdot (s_ns_{n-1})s_n.
\end{aligned}$$
where $\tau\in \Omega$ satisfies $\tau s_0=s_n\tau,\ \tau s_i=s_{n-i}\tau$, for $i=2,3,...,n-1$.

\begin{proof}
  The above expressions for $x_1$ and $x_n$ are given in  \cite{lusztig1983some}. By the length formulas 4.10 (a) and 4.10 (b) we know that above expressions for $x_1$ and $x_n$ are reduced.

 Now we show that the expression is reduced for each $x_i,\ 1\le i\le n-1$. We use induction on  $i$. It is known that for  $i=1$, the expression is reduced.  Suppose the assertion is true for $i-1$.

Notice for $1\le k\ne i\le n$, we have $s_kx_i=x_is_k$. Therefore,
$$\begin{aligned}
 x_i&=s_{i-1}\cdots s_2s_1x_1s_1s_2\cdots s_{i-1}  x_{i-1}\\
&= s_{i-1}\cdots s_2s_1s_0s_1s_2\cdots s_{n-1}s_ns_{n-1}\cdots s_{i +1}s_i(s_0s_1\cdots s_{n-1}s_ns_{n-1}\cdots s_{i}s_{i-1})^{i-1}\\
&=s_1s_2\cdots s_{i-1}x_is_{i-1}\cdots s_2s_1\\
&=s_0s_1s_2\cdots s_{n-1}s_ns_{n-1}\cdots s_{i+1}s_i(s_0s_1\cdots s_{n-1}s_ns_{n-1}\cdots s_{i}s_{i-1})^{i-1}s_{i-1}\cdots s_2s_1.
\end{aligned}$$
By commutation relation $s_ks_{k+1}s_k=s_{k+1}s_ks_{k+1},\ 1\le k\le n-2$ and the induction hypothesis, it is easy to see for $1\le k\le i-1$,
$$s_0s_1\cdots s_{n-1}s_ns_{n-1}\cdots s_{i}s_{i-1}(s_0s_1\cdots s_{n-1}s_ns_{n-1}\cdots s_{i})^{i-1-k}s_k=
(s_0s_1\cdots s_{n-1}s_ns_{n-1}\cdots s_{i})^{i-k}.$$
Thus
$$ x_i=s_{i-1}\cdots s_2s_1x_1s_1s_2\cdots s_{i-1}  x_{i-1}=(s_0s_1\cdots s_{n-1}s_ns_{n-1}\cdots s_{i})^{i }.$$
According to 4.10 (a), the number of simple reflections in the expression above is equal to $l(x_i)$, this means that the expression is reduced. The proposition follows.
\end{proof}

{\bf Theorem 4.12.} Keep the notations in 4.10. Then $W$ is an extended affine Weyl group of type $\tilde C_n\ (n\ge 2)$ and $x_1,\ ...,\ x_n$ are fundamental dominant weights in $W$.
Let $J=\{2,\ 3,\ \cdots,\ n\}$, $x_J=\displaystyle \prod_{i\in J}x_i$.
 Then $x_I\in\Omega_{qr}$.
In particular,
$x=\prod_{i\in J}x_i^{a_i}\in \Omega_{qr},$ if all $  a_i's$ are  positive integers.

\begin{proof}
By Theorem 3.8 (a),  $x_J$  is not in the lowest two-sided cell $W_{(\nu)}$ of $W$. Therefore $a(x)\le (n-1)^2+1=a(\Omega_{qr})$.

   Now we show $x_J=w\cdot w_K\cdot u$, where $K=J\cup\{0\}$ and $w_K$ is the longest element in $W_K=\{s_i\,|\, i\in K\}$.  Let $y=x_n x_{n-1}\cdots x_3, \ z:=s_0s_1\cdots s_{n-1}s_n s_{n-1}\cdots s_3s_2$. Using Proposition 4.11 we get  $x_J=y{z}^2$.
According to Lemma 2.3, $\mathscr{R}(y)=\{s_3,\ s_4,\ ...,\ s_n\}$. Note that  for any $3\le i\le n$ we have  $s_iz=zs_i$. So $\mathscr{R}(yz)\supset \{s_i\,|\, i\in J\}$. In addition, $s_0z< z$. Thus $x_J=y{z}^2=w\cdot w_K\cdot u$. Since $a(w_k)=l(w_K)=a(\Omega_{qr})$, we see $w_K\in\Omega_{qr}$. The theorem then follows from $w_K\in\Omega_{qr}$ and $x_J\not\in W_{(\nu)}$.
\end{proof}

\bigskip

{\bf 4.13} Now we deal with type $\tilde F_4$.

\def\ve{\varepsilon}
 Let $E$ a the 4-dimensional  euclidean space with the standard orthogonal basis $\ve_1,\ \ve_2,\ \ve_3,\  \ve_4$, $\Phi=\{\pm(\ve_i\pm\ve_j),\ \pm \ve_k,\ \pm\frac{1}{2}(\ve_1\pm\ve_2\pm\ve_3\pm\ve_4)\,|\, 1\le i<j\le 4,\ 1\le k\le 4\}.$ As before, $W$ is the extended affine Weyl group associated to $\Phi$, which equals $W_a$ in this case, $W_0$ is the corresponding  Weyl group, $  X$ is the subgroup consisting of all translations in $W$. The  simple refections $s_0, s_1, \cdots, s_4$ in $W$ are numbered as usual.

The Coxeter graph of $W$ is as follows ( see \cite{shi2011second} ):

\begin{center}
  \begin{tikzpicture}[scale=.6]
    \draw (-1,0) node[anchor=east]  {$\tilde F_4:$};
    \draw[thick] (0 cm,0) circle (.2 cm) node [below] {$0$};
    \draw[thick] (2 cm,0) circle (.2 cm) node [below] {$1$};
    \draw[thick] (4 cm,0) circle (.2 cm) node [below] {$2$};
    \draw[thick] (6 cm,0) circle (.2 cm) node [below] {$3$};
    \draw[thick] (8 cm,0) circle (.2 cm) node [below] {$4$};
    \draw[thick] (.2 cm,0) -- +(1.6 cm,0);
    \draw[thick] (2.2 cm,0) -- +(1.6 cm,0);
    \draw[thick] (4.2 cm,0) -- +(1.6 cm,0);
     \draw[thick] (4.2 cm,-.1) -- +(1.6 cm,0);
    \draw[thick] (6.2 cm,0) -- +(1.6 cm,0);
  \end{tikzpicture}
\end{center}

Fix a simple root system as
$$\Delta=\{\frac{1}{2}(\ve_1-\ve_2-\ve_3-\ve_4),\ \ve_4,\ \ve_3-\ve_4,\ \ve_2-\ve_3 \}.$$
Then the set of positive roots is
$$  \Phi^+=\{ \ve_i\pm\ve_j,\  \ve_k,\ \frac{1}{2}(\ve_1\pm\ve_2\pm\ve_3\pm\ve_4)\,|\, 1\le i<j\le 4,\ 1\le k\le 4\},$$
and  the set of negative roots is $\Phi^-=-\Phi^+.$ Let
$$\alpha_1=\frac{1}{2}(\ve_1-\ve_2-\ve_3-\ve_4),\ \alpha_2=\ve_4 ,\ \alpha_3=\ve_3-\ve_4,\ \alpha_{4}=\ve_{2}-\ve_3.$$
Notice that for any short root $\alpha$, $\alpha^\vee=2\alpha$, and for any long root $\beta$, $\beta^\vee=\beta$.

 The fundamental dominant weights in the weight lattice $P$ are:
$$\begin{aligned}
\lambda_1&=2\alpha_1+3\alpha_{2}+ 2\alpha_{3}+\alpha_4=\ve_1,\\
\lambda_2&=3\alpha_1+6\alpha_{2}+ 4\alpha_{3}+2\alpha_4=2\ve_1-\frac{1}{2}(\ve_1-\ve_2-\ve_3-\ve_4),\\
\lambda_3&=4\alpha_1+8\alpha_{2}+ 6\alpha_{3}+3\alpha_4=2\ve_1+\ve_2+\ve_3,\\
\lambda_4&=2\alpha_1+4\alpha_{2}+ 3\alpha_{3}+2\alpha_4=\ve_1+\ve_2.
\end{aligned}$$

Recall that $x_i=t_{\lambda_i}$.

{\bf Theorem 4.14.} Keep the notations in subsection 4.13. Then $W$ is an  affine Weyl group of type $\tilde F_4$ and $x_1,\ ...,\ x_4$ are the fundamental dominant weights in $W$.
Set $x=x_1^{a_1}x_2^{a_2}x_3^{a_3}, \ a_1,a_2 \ge 1,\ a_{3}\ge 2$, then $x\in \Omega_{qr}$.

{\bf 4.15.} {\it Proof of Theorem 4.14.}  By Theorem 3.8 (a), $x$ is not in the lowest two-sided cell $ W_{(\nu)}$ of $W$. Thus $a(x)\le 16=a(\Omega_{qr})$. Let $J=\{0,\ 1,\ 2,\ 3\}$, then $W_J=<s_i\,|\, i\in J>$ is a Weyl group of type $B_4$, whose longest element  $w_J$  is of length $16=4^2=a(\Omega_{qr})$. It suffices to show
$$x_1  x_{ 2}x_{3}^2=w\cdot w_J\cdot u.$$

By Lemma 2.3 (1), $\mathscr{R}(x_1x_{ 2}x_{3})=\{s_1,\ s_2,\  s_{3}\}$,
$$x_1x_2x_3=v\cdot w_K,$$
where $w_K$ is the longest element of the group $W_K=\langle s_1,s_2,s_3\rangle$. We only need to prove
$$w_Kx_{3}=w_J\cdot u.$$
Since $\{s_1,\ s_2,\ s_{3}\}\subset \mathscr{L}(w_Ix_{3})$, it suffices to show $s_0w_Kx_{3}\le w_Ix_{3}$, or equivalently, $x_{3}^{-1}w_Ks_0\le x_{3}^{-1}w_K$.

We have

 (a) $l(x_{3}^{-1})=42$,

  (b) $l(x^{-1}_{3}w_K)=51$.

In fact, using the length formula (1) of Iwahori-Matsumoto, we get
$$l(x_{3}^{-1})=\sum_{\alpha\in\Phi^+}|\langle -\lambda_{3},\ \alpha^\vee\rangle  |=
\sum_{\alpha\in\Phi^+}|\langle \lambda_{3},\ \alpha^\vee\rangle  |.$$
By direct computations, we see
$$l(x_3^{-1})=l(x_3)=42,\quad l(x_3^{-1}w_K)=42+9=51.$$

\medskip

Now we  show $x_{3}^{-1}w_Is_0\le x_{3}^{-1}w_K$. We need to prove $l(x_{3}^{-1}w_Ks_0)=50$.

By definition, $s_0=t_{-\alpha_0}s_{\alpha_0}$, where $-\alpha_0=\ve_1 =\lambda_1$ is the highest short root. Then
$$x_{3}^{-1}w_Ks_0=x_{3}^{-1}w_Kt_{\ve_1}s_{\alpha_0}=x_{3}^{-1}t_{w_K(\ve_1)}w_Ks_{\alpha_0}.$$

We  compute $w_K(\ve_1)$.  Since $s_i(\ve_1)=\ve_1$ if $i=2,3,4$ and  $w_K=s_1s_2s_3s_2s_1s_2s_3s_2s_3$, we have
$$\begin{aligned}
  w_K(\ve_1)&=s_1s_2s_3s_2s_1(\ve_1)=s_1s_2s_3s_2(\ve_1-\alpha_1)\\
&=s_{1}s_{ 2}s_{ 3} (\ve_1-\alpha_1-\alpha_2)=s_{1}s_{ 2}  (\ve_1-\alpha_1-\alpha_2-\alpha_3)\\
&=s_{ 1}(\ve_1-\alpha_1-2\alpha_2- \alpha_{3})\\
&= \ve_1-2\alpha_1-2\alpha_2- \alpha_{3}\\
&=\ve_1-(\ve_1-\ve_2-\ve_3-\ve_4)-2\ve_4-(\ve_3-\ve_4)\\
&= \ve_2.
\end{aligned}$$

By the length formula (1) of Iwahori-Matsumoto, we get
$$\begin{aligned}
 &\ \qquad  l(x_{3}^{-1}w_Ks_0)=l(x_{3}^{-1}t_{\ve_{2} }w_Ks_{\alpha_0})\\
&=\sum_{\substack{\alpha\in{\Phi^+}\\ s_{\alpha_0}w_K(\alpha)\in{\Phi^+}}}|\langle -\lambda_{3}+\ve_{2},\ \alpha^\vee\rangle  |+
\sum_{\substack{\alpha\in{\Phi^+}\\ s_{\alpha_0}w_K(\alpha)\in{\Phi^-}}}|\langle -\lambda_{3}+\ve_{2},\ \alpha^\vee\rangle  -1|\\
&=\sum_{\substack{\alpha\in{\Phi^+}\\ s_{\alpha_0}w_K(\alpha)\in{\Phi^+}}}|\langle -2\ve_{1}-\ve_{3},\ \alpha^\vee\rangle  |+
\sum_{\substack{\alpha\in{\Phi^+}\\ s_{\alpha_0}w_K(\alpha)\in{\Phi^-}}}|\langle -2\ve_{1}-\ve_{3} ,\ \alpha^\vee\rangle  -1|.
\end{aligned}$$

We need to figure out the positive roots mapped by $s_{\alpha_0}w_K$ into $\Phi^+$ or $\Phi^-$ respectively in detail.

We write $\alpha>0$ if $\alpha$ is  positive root, and write $\alpha<0$ when $\alpha$ is   negative root.
Since $w_K(\alpha_i)=-\alpha_i$ if $i=2,\ 3,\ 4$,
we see
$$s_{\alpha_0}(w_K):\ \ve_1\mapsto \ve_2,\ \ve_2\mapsto-\ve_1,\ \ve_3\mapsto -\ve_3,\ \ve_4\mapsto-\ve_4.$$
It is easy to get
$$\begin{aligned}
 &s_{\alpha_0}w_K(\frac{1}{2}(\ve_1-\ve_2\pm\ve_3\pm\ve_4))=\frac{1}{2}(\ve_1+\ve_2\pm\ve_3\pm\ve_4)>0,\\
&s_{\alpha_0}w_K(\ve_1)=\ve_2>0,\\
&s_{\alpha_0}w_K(\ve_k)=-\ve_k>0,\ k=3,\ 4,\\
&s_{\alpha_0}w_K(\ve_1\pm\ve_j)=\ve_2\mp\ve_j>0,\ j=3,\ 4,\\
&s_{\alpha_0}w_K(\ve_1-\ve_2)=\ve_1+\ve_2>0,\ \\
&s_{\alpha_0}w_K(\frac12(\ve_1+\ve_2\pm\ve_3\pm\ve_4))=\frac{1}{2}(-\ve_1+\ve_2\pm\ve_3\pm\ve_4)<0,\\
&s_{\alpha_0}w_K(\ve_2)=-\ve_1<0,\\
&s_{\alpha_0}w_K(\ve_2\pm\ve_j)=-\ve_1\mp\ve_j<0,\ \ j=3,\ 4,\\
&s_{\alpha_0}w_K(\ve_1-\ve_2)=\ve_1+\ve_2>0,\ \\
&s_{\alpha_0}w_K(\ve_3\pm\ve_4)=-\ve_3\mp\ve_4<0.
\end{aligned}$$

Now we compute $|\langle -2\ve_{1}-\ve_{3},\ \alpha^\vee\rangle |$and
$|\langle -2\ve_{1}-\ve_{3},\ \alpha^\vee\rangle  -1|$ in accordance with
$s_{\alpha_0}w_K(\alpha)>0$ or $<0$ respectively.

Case 1: $ \alpha\in\Phi^+,\ s_{\alpha_0}w_K(\alpha)>0$.

(1) When $\alpha=\frac12(\ve_1-\ve_2+\ve_3\pm\ve_4)$, we have
$$|\langle -2\ve_{1}-\ve_{3},\ \alpha^\vee\rangle  |=3$$
Thus, the two positive roots contribute 6 to the sum in the length formula.

(2) When $\alpha=\frac12(\ve_1-\ve_2-\ve_3\pm\ve_4)$, we have
$$|\langle -2\ve_{1}-\ve_{3},\ \alpha^\vee\rangle  |=1$$
Thus, the  two positive roots contribute 2 to the sum in the length formula.

(3) When $\alpha=\ve_1$, we have
$$|\langle -2\ve_{1}-\ve_{3},\ \alpha^\vee\rangle  |=4 .$$
Thus, this positive root contributes 4 to the sum in the length formula.

(4) When $\alpha= \ve_1+\ve_3$, we have
$$|\langle -2\ve_{1}-\ve_{3},\ \alpha^\vee\rangle  |=3. $$
Thus, this positive root contributes 3 to the sum in the length formula.

(5) When $\alpha= \ve_1+\ve_4,\ \ve_1-\ve_2,\ \ve_1-\ve_4$, we have
$$|\langle -2\ve_{1}-\ve_{3},\ \alpha^\vee\rangle  |=2. $$
Thus, the  three positive roots contribute 6 to the sum in the length formula.

(6) When $\alpha= \ve_1-\ve_3$, we have
$$|\langle -2\ve_{1}-\ve_{3},\ \alpha^\vee\rangle  |=1. $$
Thus, this positive root contributes 1 to the sum in the length formula.

\medskip

Case 2: $\alpha\in\Phi^+,\ s_{\alpha_0}w_K(\alpha)<0$.

(1) When $\alpha=\frac12(\ve_1+\ve_2+\ve_3\pm\ve_4)$, we have
$$|\langle -2\ve_{1}-\ve_{3},\ \alpha^\vee\rangle  -1|=4.$$
Thus, the  two positive roots contribute 8 to the sum in the length formula.

(2) When $\alpha=\frac12(\ve_1+\ve_2-\ve_3\pm\ve_4)$, we have
$$|\langle -2\ve_{1}-\ve_{3},\ \alpha^\vee\rangle  -1|=2.$$
Thus, the  two positive roots contribute 4 to the sum in the length formula.

(3) When $\alpha=\ve_2,\   \ve_4$, we have
$$|\langle -2\ve_{1}-\ve_{3},\ \alpha^\vee\rangle  -1|=1.$$
Thus, the  two positive roots contribute 2 to the sum in the length formula.

(4) When $\alpha=  \ve_3$, we have
$$|\langle -2\ve_{1}-\ve_{3},\ \alpha^\vee\rangle  -1|=3.$$
Thus,  this positive root contributes 3 to the sum in the length formula.

(5) When $\alpha= \ve_1+\ve_2$, we have
$$|\langle -2\ve_{1}-\ve_{3},\ \alpha^\vee\rangle  -1|=3. $$
Thus, this positive root contributes 3 to the sum in the length formula.

(6) When $\alpha= \ve_2+\ve_3,\ \ve_3\pm\ve_4$, we have
$$|\langle -2\ve_{1}-\ve_{3},\ \alpha^\vee\rangle  -1|=2. $$
Thus, the  three positive roots contribute 6 to the sum in the length formula.

(7) When $\alpha= \ve_2\pm\ve_4$, we have
$$|\langle -2\ve_{1}-\ve_{3},\ \alpha^\vee\rangle  -1|=1. $$
Thus, the  two positive roots contribute 2 to the sum in the length formula.

(8) When $\alpha= \ve_2-\ve_3$, we have
$$|\langle -2\ve_{1}-\ve_{3},\ \alpha^\vee\rangle  -1|=0. $$
Thus, this positive root contributes 0 to the sum in the length formula.

Adding up all the numbers in Case 1 and Case 2 above, finally we get:
$$\begin{aligned}
l(x_{3}^{-1}w_Ks_0)=&6+2+4+3+6+1\\
&+8+4+3+2+3+6+2+0\\
=&50=l(x_{3}^{-1}w_K)-1,
\end{aligned}$$
which is exactly what we need. Then Theorem 4.14 is proved.\qed

\bigskip
{\bf 4.16.} Now we consider type $\tilde G_2$.

 Let $E$ be the 3-dimensional  euclidean space with the standard orthogonal basis $\ve_1,\ \ve_2,\ \ve_3$,  $\Phi=\pm\{(\ve_1\pm\ve_2),\ (\ve_1\pm\ve_3),\ (\ve_2\pm\ve_3),\ 2\ve_1-\ve_2-\ve_3,\ 2\ve_2-\ve_1-\ve_3,
  2\ve_3-\ve_1-\ve_2\}.$ Choose $\Delta=\{\ve_1-\ve_2,\ -2\ve_1+\ve_2+\ve_3\}$ to be the set of simple roots. In this case, the weight lattice and the root lattice are equal, so the extended affine Weyl group $W$ equals the affine Weyl group $W_a$. We number the simple reflections $s_1,\ s_2,\ s_0$ as usual. Let $x_1,\ x_2$ be the corresponding fundamental dominant weights in $W$.

The Coxeter graph of $W$ is as follows(\cite{shi2011second}):

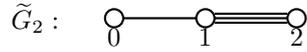
\begin{center}
  \begin{tikzpicture}[scale=.6]
    \draw (-1,0) node[anchor=east]  {$\widetilde G_2:$};
    \draw[thick] (0 cm,0) circle (.2 cm) node [below] {$0$};
    \draw[thick] (2 cm,0) circle (.2 cm) node [below] {$1$};
    \draw[thick] (4 cm,0) circle (.2 cm) node [below] {$2$};
    \draw[thick] (.2 cm,0) -- +(1.6 cm,0);
    \draw[thick] (2.2 cm,0) -- +(1.6 cm,0);
    \draw[thick] (2.2 cm,.1) -- +(1.6 cm,0);
     \draw[thick] (2.2 cm,-.1) -- +(1.6 cm,0);
  \end{tikzpicture}
\end{center}

{\bf Theorem 4.17.} Keep the notations in 4.16. So $W$ is an affine Weyl group of type $\tilde G_2$ and $x_1,\ x_2$ are fundamental dominant weights in $W$. Then a translation $x\in X^+$ is in $\Omega_{qr}$ if and only if $x={x_1}^m$ with any interger $m\ge 2$.

\begin{proof}
By Theorem 3.8 (a), $x\in\Omega_{qr}$ if and only if either $x={x_1}^m$ or $x={x_2}^m$ with some positive integer $m$.
In Example 3.6, ${x_2}^m\underset{L}\sim x_2$ for any positive integer $m$, and $x_2=s_0s_1s_2s_1s_0s_2s_1s_2s_1s_2\underset{L}\sim s_0s_1s_2s_1s_0s_2\underset{R}\sim s_0s_2$ by star actions, thus $a({x_2}^m)=2<a(\Omega_{qr}), {x_2}^m\not\in\Omega_{qr}$. In addition, $x_1\not\underset{L}\sim {x_1}^m$ while ${x_1}^m\underset{L}\sim {x_1}^2$ with $m\ge 2$ . Thus, it suffices to show  ${x_1}^2$ is in $\Omega_{qr}$  in the following.

Note that $s_1s_0s_1\in \Omega_{qr},\ a(\Omega_{qr})=3$.
Since $x_1=s_0s_1s_2s_1s_2s_1$, we see $x_1^2\in\Omega_{qr}$. Therefore, if  $m\ge 2$, then $x_1^{m}\in\Omega_{qr}$.
\end{proof}

\section{A conjecture and some related results}

{\bf 5.1.} In \cite{shi2011second}, Shi proposed a conjecture saying that the number $n_{qr}$  of left cells in the second lowest two-sided cell $\Omega_{qr}$ is half of the cardinality of the Weyl group $W_0$ and  he proved $n_{qr}\le |W_0|/2$.  He also verified his conjecture for the following cases(loc.cit.): type  $\tilde A_n$ and affine Weyl groups of rank $\le 4$. Jinhao Guo verified the conjecture for type $\tilde B_5$(see \cite{guo2019kuozhang}).

In this section we state a conjecture and provide some evidences. The conjecture could be regarded as a refinement of Shi's conjecture on the number of left cells in the second lowest two-sided cell $\Omega_{qr}$.

{\bf Conjecture 5.2.} Keep the notations in section 1. Then  $W$ is an  extended affine Weyl group associated to irreducible root system $\Phi$, $W_0$ is the corresponding Weyl group, $ X^+=\{t_\lambda\,|\,\lambda\in P^+\}\subset W$.  $x\in X^+$.  If $x$ is in the second lowest two-sided cell $\Omega_{qr}$ of $W$, then

(a) for any $w\in W_0$, we have $ wxw^{-1}\in \Omega_{qr}$;

(b) if $w,\ u\in W_0$, and $wxw^{-1}\ne uxu^{-1}$, then $wxw^{-1}\underset {R}{\not\sim} uxu^{-1}$;

(c) for any right cell $\Gamma$ in $\Omega_{qr}$ , there must  exist  some $w\in W_0$ such that $wxw^{-1}\in\Gamma$. In particular, $\Omega_{qr}$ contains $\frac{|W_0|}{2}$ right cells.

If the conjecture is true, using the results in section 4  we then obtain some representatives for each left cell in $\Omega_{qr}$. We will show that a weak version of Conjecture 5.2 (a) is true  and Conjecture 5.2 is true for type $\tilde A_{n}$ and type $\tilde G_2$.

{\bf Lemma 5.3.}  Keep the notations in section 1. Then  $W$ is  the extended affine Weyl group of the irreducible root system $\Phi$ and  $ X^+$ is the set of dominant weights in $W$. Assume that $x\in X^+$ is in the  the second lowest two-sided cell $\Omega_{qr}$. Then for any $w\in W_0$ such that $\mathscr{R}(wx)=\mathscr{R}(x)$, we have $ wx\in\Omega_{qr}$.

\begin{proof}
By Proposition 3.10 (c), there exists $\alpha_i\in{\Phi}^{+}$ such that $k(x, \alpha_i)=\langle x,\ {\alpha_i}^\vee\rangle=0$. By Proposition 3.10 (a),  $k(x^{-1}w^{-1}, \alpha_i)=\langle x^{-1}, {\alpha_i}^\vee\rangle+k(w^{-1}, \alpha_i)=-\langle x, {\alpha_i}^\vee\rangle+k(w^{-1}, \alpha_i)=k(w^{-1}, \alpha_i)$. By the length formula (1) of Iwahori-Matsumoto, we get
$$\begin{aligned}
l(xs_i)&=|\langle x,\ {\alpha_i}^\vee\rangle-1|+\sum\limits_{\substack{\alpha\in{\Phi}^+\\ \alpha\not=\alpha_i}}|\langle x,\ {\alpha}^\vee\rangle|\\
&=1+\sum\limits_{\substack{\alpha\in{\Phi}^+\\ \alpha\not=\alpha_i}}|\langle x,\ {\alpha}^\vee\rangle|\\
&= 1+l(x).
\end{aligned}$$
This implies $xs_i=s_ix$. If $k(w^{-1}, \alpha_i)\ne 0$, by Proposition 3.10 (a), $w(\alpha_i)\in{\Phi}^-$, i.e. $ws_i<w$. Then $s_i\in\mathscr{R}(wx)$, which contradicts $\mathscr{R}(wx)=\mathscr{R}(x)$. So we must have $w(\alpha_i)\in\Phi^+$. By Proposition 3.10 (a) and (c), we get  $k(x^{-1}w^{-1},  \alpha_i)=0$. This implies that $x^{-1}w^{-1}$ is not in $W_{(\nu)}$, which is equivalent to  $wx\not\in W_{(\nu)}$.

Since $wx=w\cdot x,\ x\in\Omega_{qr}$, so $a(wx)\geq a(x)$.
Consequently, $wx\in\Omega_{qr},\ a(wx)= a(x)$.
The lemma follows.
\end{proof}

The following result shows that  a weak version of Conjecture 5.2 (a) is true.

{\bf Theorem 5.4.}
Let $W$ be the extended affine Weyl group of an irreducible root system $\Phi$ of rank $n$, $ W_0\subset W$ be the  Weyl group of $\Phi$, $X$ be the set of translations (isomorphic to weight lattice), $X^+$ be the set of dominant weights, and $x_1,\ x_2,\ \cdots,\ x_n$ be the fundamental dominant weights. Set
$I_k=\{1,\ 2,\ \cdots,\ n\}\setminus\{k\}$. Assume that
$x=\prod_{i\in I_k}{x_i}^{a_i}\in X^+\cap \Omega_{qr}$. Then for any $w\in W_0$, we have
$$wxx_{I_k}w^{-1}\in \Omega_{qr},$$
where $x_{I_k}=\prod_{i\in {I_k}}x_i$.

\begin{proof}
By Theorem 3.8 (a), $wxx_{I_k}w^{-1}\not\in W_{(\nu)}$. Note that for $s=s_k$ we have $sxs=x$ and $sx_{I_k}s=x_{I_k}$.  So it is no harm to assume $ws\ge w$. By the length formula (1) of Iwahori-Matsumoto,
$$x_{I_K}=(x_{I_K}w^{-1})\cdot w.$$
Thus $wxx_{I_K}w^{-1}=(wx)\cdot( x_{I_K}w^{-1})\underset{R}{\leq} wx$.
According to Lemma 5.3, $a(wx{x_{I_K}}w^{-1})\geq a(wx)=a(x)$.
This forces that  $a(wx{x_{I_K}}w^{-1})= a(wx)$ and $wxx_{I_K}w^{-1}\in \Omega_{qr}$.
\end{proof}

In the following, we will prove Conjecture 5.2 is true for an (extended) affine Weyl group of $\tilde G_2$ or $\tilde A_n$.

{\bf Proposition 5.5.}
(\cite{Xi1994Representations}, 11.2)
Let $W_a$ be an affine Weyl group of type $\tilde G_2$ generated by $s_0, s_1, s_2$, where $s_0s_2=s_2s_0$. Then the second lowest two-sided cell of $W$ is
$$ \Omega_{qr}=\{u(i,\ j,\ k) |1\leq i, j\leq 6, k\in\mathbb{N}\},$$
where $u(i,\ j,\ k)=u_is_0s_1s_0(s_2s_1s_0)^{k}{u_j}^{-1}$, $u_1=e,\  u_2=s_2,\ u_3=s_1s_2,\ u_4=s_2s_1s_2,\
u_5=s_1s_2s_1s_2,\ u_6=s_0s_1s_2s_1s_2$. Moreover,

(a) $\mathcal{D}\cap \Omega_{qr}=\{u(i,\ i,\ 0) |1\leq i\leq 6\}$ is the set of distinguished involutions in $\Omega_{qr}$,

(b)\ $u(i,\ j,\ k)\underset{L}{\sim}u(m,\ n,\ k')$ if and only if $j=n,$

\quad \ \ $u(i,\ j,\ k)\underset{R}{\sim}u(m,\ n,\ k')$ if and only if $i=m$.

\medskip

{\bf Theorem 5.6.} Conjecture 5.2 is true for  an  affine Weyl group of $\tilde G_2$.

\begin{proof}
Keep the notations in Proposition 5.5. For simplicity, in the proof we shall write $i_1i_2\cdots i_m$ for reduced expression $s_{i_1}s_{i_2}\cdots s_{i_m}$.

Let $x$ be a dominant weight in the second lowest two-sided cell $\Omega_{qr}$,  then by Theorem 4.17, $x={x_1}^a=(012121)^a$ for some integer $a\ge 2$. Since $s_2x=xs_2$, we may assume $ws_2> w,\ w\in W_0$ when considering $wxw^{-1}$, i.e. $w\in\{e, 1, 21, 121, 2121, 12121\}$, and we denote these  six elements as $w_i,\ 1\leq i\leq 6$ correspondingly.
We have

(1)
\begin{align*}
w_1x{w_1}^{-1}&=x\\
&=012121012121(012121)^{a-2}\\
&=(01212)010(2121)(012121)^{a-2}\\
&=u_6(010)(210)^{2(a-2)}{u_5}^{-1}=u(6,\ 5,\ 2(a-2)).
\end{align*}

(2) \begin{align*}
w_2x{w_2}^{-1}&=1012121(012121)^{a-2}(01212)\\
&=(010)(210)^{2(a-1)}212\\
&=u_1(010)(210)^{2(a-1)}{u_4}^{-1}=u(1,\ 4,\ 2(a-1))
\end{align*}

(3) \begin{align*}
w_3x{w_3}^{-1}&=21012121(012121)^{a-2}(0121)\\
&=2(010)(210)^{2(a-1)}21\\
&=u_2(010)(210)^{2(a-1)}{u_3}^{-1}=u(2,\ 3,\ 2(a-1))
\end{align*}

(4)\begin{align*}
w_4x{w_4}^{-1}&=121012121(012121)^{a-2}(012)\\
&=12(010)(210)^{2(a-1)}2\\
&=u_3(010)(210)^{2(a-1)}{u_2}^{-1}=u(3,\ 2,\ 2(a-1))
\end{align*}

(5)\begin{align*}
w_5x{w_5}^{-1}&=2121012121(012121)^{a-2}(01)\\
&=212(010)(210)^{2(a-1)}e\\
&=u_4(010)(210)^{2(a-1)}{u_1}^{-1}=u(4,\ 1,\ 2(a-1))
\end{align*}

(6)\begin{align*}
w_6x{w_6}^{-1}&=12121012121(012121)^{a-2}(0)\\
&=1212(010)(210)^{2(a-2)}21210\\
&=u_5(010)(210)^{2(a-2)}{u_6}^{-1}=u(5,\ 6,\ 2(a-2))
\end{align*}

By Proposition 5.5, one can get for any $w\in W_0$, $wxw^{-1}\in\Omega_{qr}$, and  when $ws_2> w$ and $us_2> u$,  $w\not=u$, we have $wxw^{-1}\underset{R}{\not\sim}uxu^{-1}$. Part (a) and (b) of Conjecture 5.2 are proved.  Part (c) of  Conjecture 5.2 then follows from Proposition 5.5. \end{proof}

{\bf Theorem 5.7.}
Conjecture 5.2 is true for an extended affine Weyl group of type $\tilde A_n$.

In order to verify Conjecture 5.2 (b) is true for an extended affine Weyl group of type $\tilde A_n$, we need recall some definitions and results due to Shi\cite{shi1986kazhdan}.

{\bf Definition 5.8.} (\cite{shi1986kazhdan}, Definition 7.2.1)
Let $(W_a, S)$ be an affine Weyl group of type $\tilde A_n$, we say a set $M\subset W_a$ is left connected if, for any $x,\ y\in M$, there exists a sequence $y_0=x,\ y_1,\ \cdots,\ y_t=y$ in $M$ such that for every $j$ satisfying $1\leq j\leq t$, $y_j=s_{i_j}y_{j-1}$ for some $s_{i_j}\in S$. We say $M$ is right connected if $M^{-1}=\{ w^{-1}\,|\, w\in M\}$ is left connected.

{\bf Theorem 5.9. } (\cite{shi1986kazhdan}, Theorem 18.2.1)
Let $(W_a, S)$ be an affine Weyl group of type $\tilde A_n$, then

(a) any left cell in $W_a$ is left connected and is also a maximal left connected component in the two-sided cell containing it;

(b) any right cell in $W_a$ is right connected and is also a maximal right connected component in the two-sided cell containing it;

{\bf Lemma 5.10.} (\cite{shi1986kazhdan}, Lemma 18.2.5)
Let $(W, S)$  be an extended affine Weyl group of type $\tilde A_n$, $\Theta$ be a two-sided cell in $W$. Then

(a) any left connected set of $\Theta$ is  contained in some left cell of $W$,

(b) any right connected set of $\Theta$ is  contained in some right cell of $W$.

{\bf Lemma 5.11.}
Let $(W, S)$ be an extended affine Weyl group of type $\tilde A_n$, $W_0\subset W$ be the  Weyl group, $x=\prod_{i\in{I_k}}{x_i}^{a_i}\in X^{+}\cap\Omega_{qr}$ with $a_i\ge1,\ I_k=\{1,\ 2,\ \cdots,\ n\}\setminus \{k\}$. Then for any $w\in W_0$ with $ws_k> w$, we have $wx\in\Omega_{qr}$.

\begin{proof}
For any $w\in W_0$ satisfying $ws_k> w$ in $W_0$, we have $\mathscr{R}(wx)=\mathscr{R}(x)$. Then by Lemma 5.3, $wx\in\Omega_{qr}$.
\end{proof}

{\bf Corollary 5.12.}
Let $W,\ S,\ W_0,\ x,\ w$ be as in the Lemma 5.11. Let $w=s_{i_1}s_{i_2}\cdots s_{i_j}$ be a reduced expression, $w_l:=s_{i_l}s_{i_{l+1}}\cdots s_{i_j},\ 1\leq l\leq j$, then $wx{w_l}^{-1}\in\Omega_{qr}$. In particular, $wx\underset{R}{\sim}wx{w_l}^{-1}$.

\begin{proof}
By Theorem 3.3, $ wxw^{-1}\underset{LR}{\sim}x$. Thus $wxw^{-1}\in\Omega_{qr}$.
Since $wx=wx{w_{l}}^{-1}\cdot {w_{l}}=wxw^{-1}\cdot (w{w_{l}}^{-1})\cdot {w_{l}}$, according to 2.4(e),
$wx\underset{R}{\leq}wx{w_{l}}^{-1}\underset{R}{\leq}wxw^{-1}$, and
$a(wx)\geq a(wx{w_{l}}^{-1})\geq a(wx{w_{l-1}}^{-1})\geq\cdots\geq a(wxw^{-1})$.
Since $wx, wxw^{-1}\in\Omega_{qr}$, for any $1\leq l\leq j$, we have $ a(wx{w_l}^{-1})=a(wx),\ wx{w_l}^{-1}\in\Omega_{qr}$. Hence $wx\underset{R}{\sim}wx{w_l}^{-1}$.
\end{proof}

{\bf 5.13.} Proof of Theorem 5.7:
Let $(W, S)$ be an extended affine Weyl group of type $\tilde A_n$, $x$ be a dominant weight in $\Omega_{qr}$. Part (a) of Conjecture 5.2 for type $\tilde A_n$ is a special case of Theorem 3.3.

Now we prove part (b) of Conjecture 5.2. According to Theorem 4.5, there exists $1\le k\le n$ such that $x=\prod_{i\in I_k}x_i^{a_i},\  \forall a_i\ge 1$, where $I_k=\{1,\ 2,\ ...,\ n\}\setminus\{k\}$. Assume that $w, u\in W_0$ such that $wxw^{-1}\ne uxu^{-1}$.
Since $s_kxs_k=x$, it is no harm to assume that $ws_k> w$ and $us_k> u$. By Corollary 5.12,  $wx\underset{R}{\sim}wxw^{-1}$ and $\ ux\underset{R}{\sim}uxu^{-1}$. If $wxw^{-1}\underset{R}{\sim}uxu^{-1}$, then $wx\underset{R}{\sim}ux$. According to Theorem 5.9, $wx$ and $ux$ are in the same right connected set.
As a result, there exists $v\in W$ such that $wx=uxv$. Then $u^{-1}wxv^{-1}=x$.
Recall that $X$ is a normal subgroup of $W=W_0\ltimes X$, so $u^{-1}w=v$ and it commutes with $x$. This forces $v=e$ or $s_k$. If $v=s_k$, then $w=us_k$, which contradicts the assumption of $ws_k> w,\ us_k> u$.  Thus $v$ has to be equal to $e$ and $w=u$. This contradicts the assumption  $wxw^{-1}\ne uxu^{-1}$, so we must have $wxw^{-1}\not \underset{R}{\sim}uxu^{-1}$ if $wxw^{-1}\ne uxu^{-1}$.
Conjecture 5.2 (b) is proved for type $\tilde A_n$.

By Theorem 14.4.5 in \cite{shi1986kazhdan} , $\Omega_{qr}$ has  $|W_0|/2$ right cells. There are exactly  $|W_0|/2$ elements in the set $\{w\in W_0 |  ws_k> w\}$.
By the argument for part (b) of Conjecture 5.2, we see that if $w ,u\in W_0$, $ws_k> w,\ us_k> u$, and $w\ne u$, then $wxw^{-1}\not \underset{R}{\sim}uxu^{-1}$. Hence  any right cell in $\Omega_{qr}$ contains some $wxw^{-1}$. This proves part (c) of Conjecture 5.2 for type $\tilde A_n$.

Theorem 5.7 is proved.\qed

Now we show that the Conjecture 5.2 is true for the lowest two-sided cell of an extended affine Weyl group.

{\bf Definition 5.14.} ([15], section 2) A $\Phi$-tuple $X=({X_\alpha})_{\alpha\in\Phi}$ over the set $\{+,-,\bigcirc\}$ is called a sign type of type $\Phi$
 if the set $\{{X_\alpha},{X_{ -\alpha}}\}$  is either $\{\bigcirc, \bigcirc\}$ or $\{+,-\}$ for any $\alpha\in\Phi$. We see that a sign type
 $({X_\alpha})_{\alpha\in\Phi}$ is entirely determined by the $\Phi^+$ -tuple. We shall identify $({X_\alpha})_{\alpha\in\Phi^+}$
 with $({X_\alpha})_{\alpha\in\Phi}$ and call $({X_\alpha})_{\alpha\in\Phi^+}$ a sign type.
Let $\bar{\mathscr{S}}=\bar{\mathscr{S}}(\Phi)$ be the set of all sign types of type $\Phi$.
Let
\begin{multline}\notag
\mathcal G_1=
\left\{ \begin{matrix}
+&+  &+  & \bigcirc  &-   &+  &+ &\bigcirc  &\bigcirc  &-\\
+\; +, &+\; \bigcirc, &+\; -, &+\;-, &+\; -, &\bigcirc \; +, &\bigcirc\; \bigcirc, &\bigcirc\; \bigcirc, &\bigcirc\; -, &\bigcirc\; -,\\
\end{matrix}\right.\\
\left. \begin{matrix}
+ &\bigcirc &- &\bigcirc &- & -\\
-\; +, &- \;+, & - \;+,  &- \;\bigcirc, &- \;\bigcirc, &- \;-\\
\end{matrix}\right\}.
\end{multline}

\begin{multline}\notag
\mathcal G_2=
\left\{ \begin{matrix}
\bigcirc&\bigcirc&\bigcirc&-&-&-&+&\bigcirc&\bigcirc&\bigcirc\\
\bigcirc\;\bigcirc &+\;\bigcirc&+\;\bigcirc&\bigcirc\;\bigcirc&-\;\bigcirc&-\;\bigcirc&+\;\bigcirc&\bigcirc\;-&\bigcirc\;-&+\;-\\
\bigcirc, &\bigcirc, &+,&\bigcirc, &\bigcirc, &-,&+,&\bigcirc, &-,&\bigcirc, \\
\end{matrix}\right.\\
\qquad\qquad\qquad\qquad\qquad\left. \begin{matrix}
\bigcirc&+&+&-&+&+&+&-&-\\
-\;-&\bigcirc\;-&+\;-&-\;-&-\;-&+\;-&+\;-&\bigcirc\;+&\bigcirc\;+\\
-,&-,&\bigcirc,&-,&-,&-,&+,&\bigcirc,&+,\\
\end{matrix}\right.\\
\left. \begin{matrix}
-&\bigcirc&-&-&-&+\\
-\;+&+\;+&+\;+&-\;+&-\;+&+\;+\\
\bigcirc,&+,&+,&+,&-,&+\\
\end{matrix}\right\}.
\end{multline}

For any subsystem $\Phi'$ of $\Phi$,  $\Phi'^+= \Phi^+\cap\Phi'$ is a positive subsystem of $\Phi$.

Given an indecomposable positive subsystem $\Phi'^+$ of $\Phi$ of rank 2, we say that a sign type $({X_\alpha})_{\alpha\in\Phi'^+}$ is admissible if we have one of the following cases.

(1) $\Phi'^+$ has type $A_2$, say $\Phi'^+=\{\alpha, \beta, \alpha+\beta\}$. Then
\begin{center}
$X_{\alpha+\beta}$\\
$X_\alpha \qquad X_\beta$
\end{center}
belongs to $\mathcal G_1$.

(2) $\Phi'^+$ has type $B_2$, say $\Phi'^+=\{\alpha, \beta, \alpha+\beta, 2\alpha+\beta\}$. Then
\begin{center}
$X_{\alpha+\beta}$\\
$X_\alpha \qquad\qquad X_\beta$\\
$X_{2\alpha+\beta}$
\end{center}
belongs to $\mathcal G_2$.

We say that a sign type $({X_\alpha})_{\alpha\in\Phi}$ is admisssible if for any indecomposable positive subsystem $\Phi^{'+}$ of $\Phi$ of rank 2, the sign type $({X_\alpha})_{\alpha\in\Phi^{'+}}$  is admissible.
Let $\mathscr{S}=\mathscr{S}(\Phi)$ be the set of all admissible sign types of $\bar{\mathscr{S}}$.

Define a map $\zeta: {W_a}\longrightarrow \mathscr{S}$ by sending $w=(k(w, \alpha))_{\alpha\in \Phi}$ to
${X_w}=(X_\alpha)_{\alpha\in\Phi}$ such that for any $\alpha\in\Phi$,
\begin{alignat}{3}
&k(w,\alpha)>0 \Leftrightarrow X_\alpha=+,\\
&k(w,\alpha)=0 \Leftrightarrow X_\alpha=\bigcirc,\\
&k(w,\alpha)<0 \Leftrightarrow X_\alpha=-.
\end{alignat}

{\bf Theorem 5.15.} ([15], Theorem 2.1)
$\zeta({W_a})=\mathscr{S}(\Phi)$ for any indecomposable root system $\Phi$.

 Fact: The elements in $W_a$ having the same sigh type form an equivalence class of $W$, we call it an ST-class of $W_a$.

$\bf{Theorem 5.16.}$
Let $x\in W_{(\nu)}$ be a translation. Then
 
 (a) for any $w\in W_0$, $wxw^{-1}\in W_{(\nu)}$;
 
 (b) if $w, u\in W_0$, and $w\ne u$, $wxw^{-1}\not\underset{R}\sim uxu^{-1}$;

 (c) for any right cell $\Gamma$ in $W_{(\nu)}$, there must exist some $w\in W_0$ such that $wxw^{-1}\in\Gamma$.

 \begin{proof}
 (a) follows from Theorem 3.8.

 In \cite{jian1988two} Corollary 1.2, Shi proved each left cell of $W_{(\nu)}$ is an ST-class of $W$.
 For different elements $w$ and $u$ in $W_0$, there must exist some $\alpha\in\Phi$, such that $w^{-1}(\alpha)\in\Phi^+$, and $u^{-1}(\alpha)\in\Phi^-$.  We may assume that $x$ is dominant. Then $k(wx^{-1}w^{-1}, \alpha)=\langle x^{-1}, (w^{-1}(\alpha))^\vee\rangle<0$, and $k(ux^{-1}u^{-1}, \alpha)=\langle x^{-1}, (u^{-1}(\alpha))^\vee\rangle>0$. Hence $wxw^{-1}\not\underset{L}\sim ux^{-1}u^{-1}$, equivalently, $wxw^{-1}\not\underset{R}\sim uxu^{-1}$. (b) is proved.

 In \cite{jian1988two}, Shi proved that $W_{(\nu)}$ contains exactly $|W_0|$ left( respectively, right) cells. And (b) shows that for a fixed translation $x\in W_{(\nu)}$, the $|W_0|$ elements $wxw^{-1}$ with $w\in W_0$ are in different right cells, therefore, (c) is true.
 \end{proof}

\section*{Acknowledgement}
I would like to thank Professor Nanhua Xi for his supervision on this paper.

\renewcommand{\refname}{Reference}
\bibliographystyle{plain}
\bibliography{ref}

\end{spacing}
\end{document}